\documentclass[11pt,a4paper]{article}
\usepackage{amssymb,amsmath}
\usepackage{color}
\usepackage{a4wide}
\newcommand{\adots}{\mbox{\setlength{\unitlength}{1pt}
                          \begin{picture}(8,7)\put(-4,1){.}\put(0,4){.}
                          \put(4,7){.}\end{picture}}}

\newcommand {\eq} [1] {\begin{equation}\label{#1}}
\newcommand {\en} {\end{equation}}
\newcommand {\eqn}  {\begin{eqnarray}}
\newcommand {\enn} {\end{eqnarray}} 
\newcommand {\bstar}  {\begin{eqnarray*}}
\newcommand {\estar} {\end{eqnarray*}} 

\newcommand {\B} {{\bf B}}

\newcommand {\D} {{\bf D}}

\newcommand {\G} {{\bf G}}
\newcommand {\I} {{\bf I}}

\newcommand {\M} {{\bf M}}

\newcommand {\R} {{\bf R}}

\newcommand {\V} {{\bf V}}

\newcommand {\Y} {{\bf Y}}

\newcommand {\proof} {\par{\it Proof}. \ignorespaces}
\newcommand {\eproof}
      {\space
 {\ \vbox{\hrule\hbox{\vrule height1.3ex\hskip0.8ex\vrule}\hrule}}
 \par}
\newcommand {\setR}  {{\mathbb R}}
\newcommand {\setC}  {{\mathbb C}}

\newcommand {\mat} [1] {\left[\begin{array}{#1}}
\newcommand {\rix}     {\end{array}\right]}
\newtheorem{theorem}           {Theorem}
\newtheorem{lemma}    [theorem]{Lemma}

\newtheorem{corollary}[theorem]{Corollary}
\newtheorem{proposition} [theorem] {Proposition}

\newtheorem{remark}            {Remark}

\newcommand {\diag}  {\mathop{\rm diag}\nolimits}

\newcommand {\range} {\mathop{\rm range}\nolimits}

\newcommand {\wh} {\widehat}
\newcommand {\wt} {\widetilde}

\begin{document}
\title{Invariant subspace perturbations related to defective eigenvalues of $\Delta$-Hermitian and Hamiltonian matrices}
\author{Hongguo Xu\thanks{Department of Mathematics, University of Kansas, Lawrence, KS 66045, USA.
email: {\tt feng@ku.edu}. }
}
\date{}
\maketitle
\begin{abstract} 
Structured perturbation results for  invariant subspaces of $\Delta$-Hermitian and Hamiltonian matrices are provided. The invariant subspaces under consideration are associated with the eigenvalues perturbed 
from a single defective eigenvalue. The
results show how the original eigenvectors and generalized eigenvectors are involved in composing such perturbed invariant subspaces and eigenvectors. 
\end{abstract}

\section{Introduction} Eigenvalue problem of Hamiltonian matrices plays a fundamental role in systems and control \cite{BenKM05,Bye83,Lau91,Meh91,PaiV81,Zho97,ZhoDG96}. For a
Hamiltonian matrix, its eigenvalues are symmetric about the imaginary axis and its invariant subspaces
normally have special structures as well, e.g., \cite{FreMX02,LinMX99}. These structures and symmetry are crucial for solving application problems including optimal and robust controller
constructions \cite{BenBLMX11,BenBMX07,Gah94,Lau91,Meh91,ZhoDG96},  H$_2$ and H$_\infty$ norm computations \cite{BoyB90,BoyBK89,Sch90}, and distance and robustness problems \cite{AlaBKMM11,Bye88,FazGL21,Gri04}. In these applications, eigenvalue and  invariant subspace perturbation theory plays an important role, in particular those about the imaginary eigenvalues, e.g., \cite{MehX08,MehX24}.

In \cite{MehX24}, perturbation results about imaginary eigenvalues were developed for a special type of Hamiltonian matrices with  special Hamiltonian perturbations. In this paper we will provide general perturbation results for the eigenvalues of a Hamiltonian matrix that  are defective, in the situation that the perturbed matrix is still Hamiltonian. 
The results are derived based on the perturbation theory for general matrices provided in \cite{Xu25}. Since Hamiltonian matrices can be transformed to $\Delta$-Hermitian matrices (which will be defined later), we will first derive the perturbation results for $\Delta$-Hermitian matrices. We then use the   results to establish the perturbation theory for Hamiltonian matrices. 

The paper is organized as follows. In Section~\ref{sec2}, some results from \cite{Xu25} are provided, which are needed for establishing the structured perturbation theory. Definitions, structured canonical forms, and basic properties of $\Delta$-Hermitian and Hamiltonian matrices are given as well. The perturbation results and derivations for $\Delta$-Hermitian matrices are given in Section~\ref{sec3} and the results for Hamiltonian matrices
are in Section~\ref{sec4}. Our conclusions are contained in Section~\ref{sec5}.

Throughout the paper, $\setR, \setC$ represent the real and complex fields, respectively. 
$\setC^{m\times n} $ is the space of  $m\times n$ complex 
matrices, and $\setC^n$ is the $n$-dimensional complex column vector space.
For a matrix $B$, $B^\ast$  and $B^T$ are the adjoint (conjugate transpose) and transpose of $B$, respectively.
For a square matrix $A$, $\Lambda(A)$ is the spectrum of $A$. We
denote by $I_k$ and $0_k$ the $k\times k$ identity and zero matrices, respectively. For a positive number $d$, the notation $O(t^d)$ means
a term (or matrix) that converges to $0$ at the same rate as $t^d$ and
$o(t^d)$ means a term (or matrix) that converges  faster than $t^d$, both when $t\to 0$. 
We also use notations as $O\left(\mat{ccc}t^{d_{11}}&\ldots&t^{d_{1q}}\\
\vdots&\vdots&\vdots\\t^{d_{p1}}&\ldots&t^{d_{pq}}\rix\right)$ to show the
 orders of different blocks (entries) in a matrix. $t\in \setR$ is always a positive
and sufficiently small parameter. 
For a given $r$-dimensional subspace $V\subseteq \setC^n$,  we call $X\in\setC^{n\times r}$ a
{\em  basis matrix} of $V$  if $\range X=V$. Finally, $\imath=\sqrt{-1}$.
\section{Preliminaries}\label{sec2}
 Consider the invariant subspaces of the matrix $\lambda I_p+N+tB(t)$, where
$B(t)=B+O(t)$ for some constant matrix $B$, 
\eq{defn}
N:=N(m)=\diag(N_1,N_{2},\ldots,N_{m}),
\en
and
\[
N_j=\mat{cccc}0&I_{s_j}&&\\&0&\ddots&\\&&\ddots&I_{s_j}\\&&&0\rix_{j\times j},
\quad j=1,\ldots,m
\]
with integers $s_1,\ldots,s_{m-1}\ge 0$,  $s_{m}>0$, and
\[
p=s_1+2s_2+\ldots+m s_{m}.
\]
By convention, $N_1=0_{s_1}$. This indicates that $\lambda I_p+N$ has 
$s_j$ Jordan blocks of size $j\times j$ for $j=1,\ldots,m$.
Because of the term $tB(t)$, the original eigenvalue $\lambda$ will be perturbed with magnitudes of fractional orders of $t$ that depend on the sizes of the Jordan blocks. 

We focus on the perturbations of the
eigenvalues corresponding to a particular size of Jordan blocks 
and the corresponding invariant subspaces. 

Partition
\[
B=\mat{ccc}B_{11}&\ldots&B_{1m}\\\vdots&\ddots&\vdots\\B_{m1}&\ldots&B_{mm}\rix,\quad
B_{ij}=\mat{ccc}B_{11}^{(ij)}&\ldots&B_{1j}^{(ij)}\\\vdots&\ddots&\vdots\\
B_{i1}^{(ij)}&\ldots&B_{ij}^{(ij)}\rix
\]
conformably with the block forms of $N$ and $N_i,N_j$, respectively, for $1\le i,j\le m$. 
Let $W_m=B_{m1}^{(mm)}$ and for any
$k\in\{1,\ldots,m-1\}$, define 
\eq{defw}
W_k=
\mat{c|ccc}B_{k1}^{(kk)}&B_{k1}^{(k,k+1)}&\ldots&B_{k1}^{(km)}\\\hline
B_{k+1,1}^{(k+1,k)}&B_{k+1,1}^{(k+1,k+1)}&\ldots&B_{k+1,1}^{(k+1,m)}\\
\vdots&\vdots&\ddots&\vdots\\
B_{m1}^{(mk)}&B_{m1}^{(m,k+1)}&\ldots&B_{m1}^{(mm)}\rix
=\mat{c|c}B_{k1}^{(kk)}&W_{k,k+1}\\\hline W_{k+1,k}&W_{k+1}\rix.
\en
The following perturbation results are given
in \cite{Lid66,MorBO97,Xu25}.

\begin{lemma}\label{lem1}
Consider the matrix $\lambda I_p+N+tB(t)$  with $B(t)=B+O(t)$ and $N:=N(m)$. Let $W_k$ be defined in (\ref{defw}).
Suppose $\rho\in\{1,\ldots,m\}$ and $W_{\rho+1}$ is invertible when $\rho<m$. Let
$\gamma_1^{(\rho)},\ldots,\gamma_{s_\rho}^{(\rho)}$
be the $s_\rho$ eigenvalues of the matrix
\eq{defs}
S_\rho=\left\{\begin{array}{ll}
B_{\rho1}^{(\rho\rho)}+W_{\rho,\rho+1}G_{\rho}&\rho<m\\
B_{m1}^{(mm)}&\rho=m,
\end{array}
\right.,\qquad G_\rho=-W_{\rho+1}^{-1}
W_{\rho+1,\rho}=:\mat{c}G_{\rho+1,\rho}\\\vdots\\G_{m\rho}\rix.
\en
Assume $t$ is sufficiently small.
\begin{enumerate}
\item[(a)] 
For each $i\in\{1,\ldots,s_\rho\}$, let
$
\mu_{i1}^{(\rho)},\ldots,\mu_{i\rho}^{(\rho)}
$
be  the $\rho$th roots of $\gamma_i^{(\rho)}$, i.e., 
$(\mu_{ij}^{(\rho)})^\rho=\gamma_i^{(\rho)}$, for $j=1,\ldots,\rho$. Then
$\lambda I_p+N+tB(t)$ has $\rho s_\rho$ eigenvalues 
$$\lambda_{ij}^{(\rho)}(t)=\lambda+t^{1/\rho}\mu_{ij}^{(\rho)}+o(t^{1/\rho}),\quad
i=1,\ldots, s_\rho,\quad j=1,\ldots,\rho.
$$
\item[(b)] Suppose all $W_1,\ldots,W_m$ are invertible. Let 
\eq{deftheta}
\Theta_\rho(S_\rho)=\mat{cccc}0&I_{s_\rho}&&\\&0&\ddots&\\&&\ddots&I_{s_\rho}\\S_\rho&&&0\rix_{\rho\times\rho},
\en
where $S_\rho$ is defined in (\ref{defs}).
Then there is a matrix
\eq{defx}
X_\rho(t)=\mat{ccccc}X_{1\rho}^T(t)&\ldots&X_{\rho\rho}^T(t)&\ldots&X_{m\rho}^T(t)\rix^T+O(t),
\en
\[
X_{k\rho}(t)=
\left\{\begin{array}{ll}
O(\mat{cccc}t^{1-k/\rho}&t^{1-(k-1)/\rho}&\ldots&t^{1-1/\rho}\rix^T)&k<\rho\\
\mat{cc}I_{s_\rho}&0\\0&0\rix+\diag(0_{s_\rho},t^{1/\rho}I_{s_\rho},\ldots,t^{1-1/\rho}I_{s_\rho})&
k=\rho\\
\mat{cc}G_{k\rho}&0\\0&0\rix
+O(\mat{ccccccc}t^{1/\rho}&t^{1/\rho}&t^{2/\rho}&\ldots&t&\ldots&t\rix^T)& k>\rho.
\end{array}\right.
\]
where $G_{\rho+1,\rho},\ldots,G_{m\rho}$ are given in (\ref{defs}), such that
\[
(\lambda I_p+N+tB(t))X_\rho(t)=X_\rho(t)(\lambda I_{\rho s_\rho}+t^{1/\rho}\wt \Theta_\rho(t))
\]
for some $\wt \Theta_\rho(t) =\Theta_\rho(S_\rho)+O(t^{1/\rho})$.
\end{enumerate}
\end{lemma}
\proof The  results can be derived by applying  those in \cite{Lid66,MorBO97,Xu25} to the special case when the original matrix has a single eigenvalue $\lambda$ and is already in the Jordan canonical form.
\eproof
We point out that $S_\rho$ is the Schur complement of $W_\rho$ corresponding to 
$W_{\rho+1}$.
Note $\Lambda(\Theta_\rho(S_\rho))=\{\mu^{(\rho)}_{ij}\}$.
Note also if $W_1,\ldots,W_m$ are invertible, $S_\rho$ is invertible  for any $\rho$. Then 
$\Theta_\rho(S_\rho)$ is invertible as well, and all $\mu_{ij}^{(\rho)}$,
$1\le i\le s_\rho$, $1\le j\le \rho$, are nonzero.
\begin{remark}\label{rem0}\rm The assumption that $W_1,\ldots,W_m$ are invertible
in Lemma~\ref{lem1}(b) is generic, and is necessary  for deriving 
 the invariant subspace perturbation results, e.g., \cite{Xu25}.
\end{remark}

\medskip
The following result shows the perturbation behaviors of the invariant subspaces corresponding to a subset of $\Lambda(\lambda I_{\rho s_\rho}+t^{1/\rho}\wt \Theta_\rho(t))$. Suppose $F\in\setC^{n\times \ell}$ is a full column rank matrix and 
\[
\Theta_\rho(S_\rho)F=F\Omega,
\]
where $S_\rho$ and $\Theta_\rho(S_\rho)$ are given in (\ref{defs}) and (\ref{deftheta}),
and $\Omega\in\setC^{\ell\times \ell}$. 
Clearly $\Lambda(\Omega)\subseteq\Lambda(\Theta_\rho(S_\rho))$.
Based on the block structure of $\Theta_\rho(S_\rho)$, we have the expression
\eq{deffq}
F=\mat{cccc}Q^T&(Q\Omega)^T&\ldots&(Q\Omega^{\rho-1})^T\rix^T
\en
for some matrix $Q\in\setC^{s_\rho\times \ell}$ that satisfies 
\[
S_\rho Q=Q\Omega^\rho.
\]
\begin{lemma}\label{lem2}
Assume the conditions of Lemma~\ref{lem1} hold and $
\Lambda(\Omega)\cap(\Lambda(\Theta_\rho(S_\rho))\backslash\Lambda(\Omega))=\emptyset.$ 
Then
\[
(\lambda I_p+N+tB(t))X_\Omega(t)=X_\Omega(t)(\lambda I_{\ell}+t^{1/\rho}\wt \Omega(t)),
\]
where $\wt \Omega(t) =\Omega+O(t^{1/\rho})$, 
\[
X_\Omega(t)=\mat{ccccc}X_{\Omega,1}^T(t)&\ldots&X_{\Omega,\rho}^T(t)&\ldots&X_{\Omega,m}^T(t)\rix^T+O(t),
\]
\[
X_{\Omega,k}(t)=\left\{
\begin{array}{ll}O(\mat{cccc}t^{1-k/\rho}&t^{1-(k-1)/\rho}&\ldots&t^{1-1/\rho}\rix^T)&k<\rho\\
\mat{cccc}Q^T&(Q(t^{1/\rho}\Omega))^T&\ldots&(Q(t^{1/\rho}\Omega)^{\rho-1})^T\rix^T&\\
\qquad\qquad\qquad\qquad\qquad\quad+O(\mat{cccc}t^{1/\rho}&t^{2/\rho}&\ldots&t\rix^T)&k=\rho\\
\mat{cccc}(G_{k\rho}Q)^T&0&\ldots&0\rix^T&\\
\qquad+O(\mat{cccccccc}t^{1/\rho}&t^{1/\rho}&t^{2/\rho}&\ldots&t^{1-1/\rho}&t&\ldots&t\rix^T)& k>\rho,
\end{array}\right.
\]
and $Q$ is given in (\ref{deffq}). 
\end{lemma}
\proof
See \cite{Xu25}.
\eproof

\medskip 
Let $\Delta$ be Hermitian and invertible. A matrix $C$ is {\em $\Delta$-Hermitian} if
\[
\Delta C=(\Delta C)^\ast=C^\ast \Delta.
\]

The structured Jordan canonical form of a $\Delta$-Hermitian is well known.
\begin{proposition}\label{pro1} 
Suppose $C$ is $\Delta$-Hermitian. There exists an invertible matrix $U$ such that
\bstar
U^{-1}CU&=&\diag(C_1,\ldots,C_r,C_{r+1},\ldots,C_{q})=:\wh C,\\
U^\ast \Delta U&=&\diag(\Delta_1,\ldots,\Delta_r,\Delta_{r+1},\ldots,\Delta_q)
=:\wh \Delta=\wh \Delta^\ast=\wh \Delta^{-1}.
\estar
\begin{enumerate}
\item[(i)] For each $j=1,\ldots,r$, $C_j$ contains a pair of complex conjugate eigenvalues
$\lambda_j,\bar\lambda_j$, 
\[
C_j=\mat{cc}\lambda_j I_{p_j}+N(m_j)&0\\0&(\lambda_j I_{p_j} +N(m_j))^\ast\rix,\quad
\Delta_j=\mat{cc}0&I_{p_j}\\I_{p_j}&0\rix=\Delta_j^\ast=\Delta_j^{-1},
\]
where $N(m_j)$ is defined in (\ref{defn}) with $s^{(j)}_k$ size $k\times k$ Jordan blocks
for $k=1,\ldots,m_j$, and $p_j=\sum_{k=1}^{m_j}ks^{(j)}_{k}$.

Moreover, $\lambda_{1},\ldots,\lambda_r$ are distinct.
\item[(ii)] For each $j=r+1,\ldots,q$, $C_j$ contains a real eigenvalue $\alpha_j$, 
\[
C_j=\alpha_j I_{p_j}+N(m_j),\quad \Delta_j=\diag(\Gamma_1^{(j)},\Gamma_2^{(j)},\ldots,
\Gamma_{m_j}^{(j)})=\Delta^\ast_j=\Delta_j^{-1},
\]
and
\[
\Gamma_1^{(j)}=\Sigma_1^{(j)},\quad
\Gamma_2^{(j)}=\mat{cc}0&\Sigma_2^{(j)}\\\Sigma_2^{(j)}&0\rix_{2\times 2},\quad\ldots,\quad
\Gamma_{m_j}^{(j)}=\mat{ccc}&&\Sigma_{m_j}^{(j)}\\&\adots&\\\Sigma_{m_j}^{(j)}&&\rix_{m_j\times m_j}.
\]
\[
\Sigma_k^{(j)}=\mat{cc}I_{t^{(j)}_k}&0\\0&-I_{s^{(j)}_k-t^{(j)}_k}\rix,\qquad k=1,\ldots,m_j,
\]
where $N(m_j)$ and $p_j$ are defined in the same way as in (i) and $t_k^{(j)}\ge 0$ for 
$k=1,\ldots,m_j$. The diagonal entries of $\Sigma_k^{(j)}$ 
(which are $\pm 1$) are the {\em sign signatures} of $C$
corresponding to the $k\times k$ Jordan blocks of $\alpha_j$.

Moreover, $\alpha_{r+1},\ldots,\alpha_q$ are distinct.
\end{enumerate}

\end{proposition}
\proof This is a slightly modified form of those in \cite{GohLR83,MehX99,Tho91}.
\eproof
\medskip
Define $J_n=\mat{cc}0&I_n\\-I_n&0\rix$. A matrix $H\in\setC^{2n\times 2n}$ is 
{\em Hamiltonian} if $J_nH=(J_nH)^\ast$. 
\begin{proposition}\label{pro2} 
Suppose $H$ is Hamiltonian. There exists an invertible matrix $U$ such that
\[
U^{-1}HU=\diag(H_1,\ldots,H_r,H_{r+1},\ldots,H_{q}),\quad
U^\ast J_n U=\diag(J_{p_1},\ldots,J_{p_r},\imath\Delta_{r+1},\ldots,\imath\Delta_q).
\]
\begin{enumerate}
\item[(i)] For each $j=1,\ldots,r$, $H_j$ contains a pair of nonimaginary eigenvalues
$\lambda_j,-\bar\lambda_j$, 
\[
H_j=\mat{cc}\lambda_j I_{p_j}+N(m_j)&0\\0&-(\lambda_j I_{p_j} +N(m_j))^\ast\rix,
\]
where $N(m_j)$ and $p_j$ are as defined in Proposition~\ref{pro1}(i).

Moreover, $\lambda_{1},\ldots,\lambda_r$ are distinct.
\item[(ii)] For each $j=r+1,\ldots,q$, $H_j$ contains an imaginary eigenvalue $\imath\alpha_j$, 
\[
H_j=\imath(\alpha_j I_{p_j}+N(m_j)),\quad 
\Delta_j=\diag(\Gamma_1^{(j)},\Gamma_2^{(j)},\ldots,\Gamma^{(j)}_{m_j}),
\]
where $N(m_j)$ and $\Delta_j$ are the same as that in
Proposition~\ref{pro1}(ii). 
The diagonal entries of $\Sigma_k^{(j)}$ in $\Gamma_k^{(j)}$ are the {\em structure inertia
indices} of $H$ corresponding to the $k\times k$ Jordan blocks of $\imath \alpha_j$.

Moreover, $\imath\alpha_{r+1},\ldots,\imath\alpha_q$ are distinct.
\end{enumerate}
\end{proposition}
\proof
This is a slightly modified form of those in  \cite{FreMX02,LanR95,LinMX99}.
\eproof
\section{Invariant subspace perturbations of $\Delta$-Hermitian matrices}~\label{sec3}
We investigate how the invariant subspaces change when a $\Delta$-Hermitian matrix
$C$ is perturbed to $C+tD$, where $D$ is another given $\Delta$-Hermitian matrix
and $t$ is sufficiently small.
For an invariant subspace corresponding to certain eigenvalues of $C$ that are different from the rest, the perturbation results are well known, e.g.,
\cite{SteS90}. Our focus is about the invariant subspaces of $C+tD$
corresponding to the eigenvalues that are perturbed from a single defective eigenvalue of $C$
with the same fractional order.  We consider 
\begin{enumerate}
\item[(a)] a single defective nonreal eigenvalue $\lambda$ (together with $\bar\lambda$ due to the pairing)
\item[(b)] a single defective real eigenvalue $\alpha$.
\end{enumerate}

From Proposition~\ref{pro1}, 
there is an invertible matrix $U=\mat{cc}U_T&U_C\rix$ such that
\eq{strd}
CU=U\mat{cc}C_T&0\\0&C_C\rix,\quad U^\ast \Delta U=\mat{cc}\Delta_T&0\\0&\Delta_C\rix.
\en
In case (a),
\eq{defcda}
C_T=\mat{cc}\lambda I_p+N&0\\0&(\lambda I_p+N)^\ast\rix,\quad \Delta_T=\mat{cc}0&I_p\\I_p&0\rix,
\en
where $\lambda$ is nonreal and $\lambda,\bar\lambda\notin\Lambda(C_C)$, 
$N:=N(m)$ is defined in (\ref{defn}) and $p=\sum_{k=1}^mks_k$. 
In case (b),
\eq{defcdb}
C_T=\alpha I_p+N,\quad \Delta_T=\diag(\Gamma_1,\ldots,\Gamma_m),
\en
where $\alpha\in\setR$ and $\alpha\notin\Lambda(C_C)$,  $N:=N(m)$ is defined in (\ref{defn})
with $p=\sum_{k=1}^mks_k$, and $\Gamma_1,\ldots,\Gamma_m$ are defined as those 
$\Gamma_k^{(j)}$ in Proposition~\ref{pro1} with $t^{(j)}_k$, $s^{(j)}_k$ replaced by
$t_k$ and $s_k$. In this case, $\Delta_T$ is unitarily similar  to $\diag(I_{\tilde s},-I_{\tilde t})$ with
\eq{defst}
\tilde s=\left\{\begin{array}{ll} \sum_{j=0}^{q-1}(t_{2j+1}+js_{2j+1})+\sum_{j=1}^{q} js_{2j}&m=2q\\ \sum_{j=0}^{q}(t_{2j+1}+js_{2j+1})+\sum_{j=1}^{q} js_{2j}&m=2q+1\end{array}
\right.\qquad \tilde t=p-\tilde s.
\en
From the second relation in (\ref{strd}),
\[
U^{-1}=\mat{cc}\Delta_T&0\\0&\Delta_C\rix^{-1}U^\ast \Delta.
\]
Define
\eq{defd}
\mat{cc}D_{11}&D_{12}\\D_{21}&D_{22}\rix:=U^{-1}DU
=\mat{cc}\Delta_T&0\\0&\Delta_C\rix^{-1}U^\ast \Delta DU.
\en
One has
\[
U^{-1}(C+tD)U=\mat{cc}C_T+tD_{11}&tD_{12}\\tD_{21}&C_C+tD_{22}\rix,
\]
which is $\diag(\Delta_T,\Delta_C)$-Hermitian.

Consider the Riccati equation
\eq{ric}
(C_C+tD_{22})Y(t)-Y(t)(C_T+tD_{11})-tY(t)D_{12}Y(t)+tD_{21}=0.
\en
Since $\Lambda(C_C)\cap\Lambda(C_T)=\emptyset$, when $t$ is sufficiently small it has a solution $Y(t)
=O(t).$
Note in both cases, $\Delta_T^{-1}=\Delta_T$.
By taking the adjoint on both sides
of (\ref{ric}) and using the fact that $U^{-1}(C+tD)U$ is $\diag(\Delta_T,\Delta_C)$-Hermitian,
one has
\[
(C_T+tD_{11})\Delta_TY^\ast(t)\Delta_C-\Delta_TY^\ast(t)\Delta_C(C_C+tD_{22})-tD_{12}
+t\Delta_TY^\ast(t)\Delta_CD_{21}\Delta_TY^\ast(t)\Delta_C=0.
\]
Let
\[
\wh Z(t)=\mat{cc}I&-\Delta_TY^\ast(t)\Delta_C\\Y(t)&I\rix.
\]
One has 
\[
(C+tD)U\wh Z(t)
=U\wh Z(t)
\mat{cc}\wh C_T(t)&0\\0&\wh C_C(t)\rix,\quad
(U\wh Z(t))^\ast \Delta(U\wh Z(t))=\mat{cc}\wh \Delta_T(t)&0\\0&\wh \Delta_C(t)
\rix,
\]
where
\[
\wh C_T(t)=C_T+tD_{11}+tD_{12}Y(t),\qquad \wh C_C(t)=C_C+tD_{22}-tD_{21}\Delta_TY^\ast(t)\Delta_C,
\]
and by using $\Delta_T^\ast=\Delta_T=\Delta_T^{-1}$, $\Delta_C^\ast=\Delta_C$, 
\[\wh \Delta_T(t)=\Delta_T+Y^\ast(t)\Delta_CY(t),\quad
\wh\Delta_C(t)=\Delta_C+\Delta_CY(t)\Delta_TY^\ast(t)\Delta_C.
\]
It is straightforward to verify
\eq{invhz}
\wh Z^{-1}(t)=\mat{cc}\wh \Delta_T(t)&0\\
0&\wh \Delta_C(t)\rix^{-1}\mat{cc}\Delta_T&0\\0&\Delta_C\rix \mat{cc}I&\Delta_TY^\ast(t)\Delta_C
\\-Y(t)&I\rix.
\en

We now show that the matrix $\wh \Delta_T(t)$ is congruent to $\Delta_T$ with a matrix close to $I$.
For $\Delta_T$ there is a unitary matrix $Q$ such that
\[
Q^\ast \Delta_TQ=\diag(I_{\tilde s}, -I_{\tilde t}),
\]
where $\tilde s$ and $\tilde t$ are both $p$ in case (a) and are given in (\ref{defst}) in case (b). Partition
\[
Q^\ast Y^\ast(t)\Delta_CY(t)Q=\mat{cc}Y_{11}(t)&Y_{12}(t)\\Y_{12}^\ast(t)&Y_{22}(t)\rix
=O(t^2),
\]
where $Y_{11}(t)=Y_{11}^\ast (t)\in\setC^{\tilde s\times \tilde  s}$, 
$Y_{22}(t)=Y_{22}^\ast(t)\in\setC^{\tilde t\times \tilde t}$. When $t$ is sufficiently 
small, the matrix
\[
Z_{12}(t)=-(I_{\tilde s}+Y_{11}(t))^{-1}Y_{12}(t)=O(t^2)
\]
is well defined, and there are matrices $Z_{11}(t)=I_{\tilde s}+O(t^2)$ and 
$Z_{22}(t)=I_{\tilde t}+O(t^2)$ such that
\bstar
Z_{11}^\ast(t)Z_{11}(t)&=&I_{\tilde s}+Y_{11}(t)\\
Z_{22}^\ast(t)Z_{22}(t)&=&I_{\tilde t}-Y_{22}(t)+Y_{12}^\ast(t)(I_{\tilde s}+Y_{11}(t))^{-1}Y_{12}(t).
\estar
Define
\[
Z_1(t)=Q\mat{cc}I_{\tilde s}&Z_{12}(t)\\0&I_{\tilde t}\rix\mat{cc}Z_{11}^{-1}(t)&0\\0&Z_{22}^{-1}(t)\rix
Q^\ast=I+O(t^2).
\]
One has
\eq{eqdt}
Z_1^\ast(t)\wh\Delta_T(t)Z_1(t)=\Delta_T.
\en
Similarly, there is a $Z_2(t)=I+O(t^2)$ such that
\eq{eqdc}
Z_2^\ast(t)\wh\Delta_C(t)Z_2(t)=\Delta_C.
\en
Define
\[
Z(t)=\wh Z(t)\mat{cc}Z_1(t)&0\\0&Z_2(t)\rix.
\]
Then from (\ref{invhz}), (\ref{eqdt}), and (\ref{eqdc}), one has
\[
Z^{-1}(t)=\mat{cc}\Delta_TZ^\ast_1(t)\Delta_T&0\\
0&\Delta_C^{-1}Z_2^\ast(t)\Delta_C\rix
\mat{cc}I&\Delta_TY^\ast(t)\Delta_C\\-Y(t)&I\rix,
\]
and 
\eqn
\nonumber
&&
(C+tD)UZ(t)=UZ(t)\mat{cc}\wt C_T(t)&0\\0&\wt C_C(t)\rix\\
\nonumber
&&
\wt C_T(t)=Z_1^{-1}(t)\wh C_T(t)Z_1(t)=Z_1^{-1}(t)(C_T+tD_{11}+tD_{12}Y(t))Z_1(t)\\
\label{nrcase}
&&\\
\nonumber
&&
\wt C_C(t)=Z_2^{-1}(t)\wh C_C(t)Z_2(t)=Z_2^{-1}(t)(C_C+tD_{22}-tD_{21}\Delta_TY^\ast(t)\Delta_C)Z_2(t)\\
\nonumber
&&
(UZ(t))^\ast \Delta(UZ(t))=Z^\ast(t)\mat{cc}\Delta_T&0\\0&\Delta_C\rix Z(t)
=\mat{cc}\Delta_T&0\\0&\Delta_C\rix.
\enn
Obviously, $\diag(\wt C_T(t),\wt C_C(t))$ is $\diag(\Delta_T,\Delta_C)$-Hermitian.
Note that the columns of 
\[
U\mat{c}I\\Y(t)\rix Z_1(t),\qquad U\mat{c}-\Delta_TY^\ast(t)\Delta_C\\I\rix Z_2(t)
\]
span the invariant subspaces of $C+tD$ corresponding to the eigenvalues of 
$\wt C_T(t)$ and $\wt C_C(t)$, respectively. 

In the following, we
consider the perturbations of the invariant subspaces corresponding to
the eigenvalues perturbed from a single defective eigenvalue of $C$
of the same fractional order. 
\subsection{Nonreal eigenvalue case}\label{sub11}
In this case $C_T$ and $\Delta_T$ are given in (\ref{defcda}).
We need to transform $\wt C_T(t)$ further to a block diagonal form as $C_T$. 
Because $Z_1(t)=I+O(t^2)$ and $Y(t)=O(t)$, we have the expression
\[
\wt C_T(t)=C_T+t\wt D_{11}(t),\qquad \wt D_{11}(t)=D_{11}+O(t).
\]
Recall $\wt C_T(t)$ and $C_T$ are $\Delta_T$-Hermitian. So is $\wt D_{11}(t)$. 
In this case $\Delta_T=\mat{cc}0&I_p\\I_p&0\rix$. Partition
\[
D_{11}=\mat{cc} M_{11}& M_{12}\\M_{21}&M_{11}^\ast\rix,\qquad
\wt D_{11}(t)=\mat{cc}\wt M_{11}(t)&\wt M_{12}(t)\\\wt M_{21}(t)&\wt M_{11}^\ast(t)\rix,
\]
conformably, where $M_{12}=M_{12}^\ast$, $M_{21}=M_{21}^\ast$,
$\wt M_{12}(t)=\wt M_{12}^\ast(t)$, and $\wt M_{21}(t)=\wt M_{21}^\ast(t)$.
As before, we perform a structure preserving transformation on $\wt C_T(t)$ to eliminate the
(1,2) and (2,1) blocks. Since $\lambda\ne \bar\lambda$, 
when $t$ is sufficiently small, the Riccati equation
\[
(\lambda I_p+N+t\wt M_{11}(t))^\ast\wt Z_{21}(t)-\wt Z_{21}(t)(\lambda I_p+N+t\wt M_{11}(t))
-t\wt Z_{21}(t)\wt M_{12}(t)\wt Z_{21}(t)+t\wt M_{21}(t)=0
\]
has a unique skew-Hermitian solution $\wt Z_{21}(t)=-\wt Z_{21}^\ast(t)=O(t)$.
Let $\wt Z_{12}(t)=-\wt Z_{12}^\ast(t)=O(t)$ be the unique skew-Hermitian solution of the Sylvester equation
\bstar
(\lambda I_p+N+t(\wt M_{11}(t)+\wt M_{12}(t)\wt Z_{21}(t)))\wt Z_{12}(t)
-\wt Z_{12}(t)(\lambda I_p+N+t(\wt M_{11}(t)+\wt M_{12}(t)\wt Z_{21}(t)))^\ast&&\\
+t\wt M_{12}(t)
=0.&&
\estar
Then 
\[
\wt Z_1(t)=\mat{cc}I_p&0\\\wt Z_{21}(t)&I_p\rix\mat{cc}I_p&\wt Z_{12}(t)\\0&I_p\rix
=\mat{cc}I_p&\wt Z_{12}(t)\\\wt Z_{21}(t)&I_p+\wt Z_{21}(t)\wt Z_{12}(t)\rix =I_{2p}+O(t)
\]
satisfies
\eq{proz1}
\wt Z_1^\ast(t)\Delta_T\wt Z_1(t)=\Delta_T.
\en
and
\[
\wt C_T(t)\wt Z_1(t)=\wt Z_1(t)\mat{cc}\wt C_1(t)&0\\0&\wt C_1^\ast(t)\rix,
\qquad \wt C_1(t)=\lambda I_p+N+t(\wt M_{11}(t)+\wt M_{12}(t)\wt Z_{21}(t)).
\]
Hence,
\eqn
\nonumber
&&(C+tD)UZ(t)\mat{cc}\wt Z_1(t)&0\\0&I\rix
=UZ(t)\mat{cc}\wt Z_1(t)&0\\0&I\rix\mat{cc|c}\wt C_1(t)&0&0\\
0&\wt C_1^\ast(t)&0\\\hline
0&0&\wt C_C(t)\rix\\
\label{nonimag}
&&\\
&&\nonumber
\left(UZ(t)\mat{cc}\wt Z_1(t)&0\\0&I\rix\right)^\ast \Delta
UZ(t)\mat{cc}\wt Z_1(t)&0\\0&I\rix=\mat{c|c}\Delta_T&0\\\hline0&\Delta_C\rix.
\enn

We now consider the invariant subspace perturbations of $\wt C_1(t)$ and $\wt C_1^\ast(t)$
corresponding to the eigenvalues perturbed from $\lambda$ and $\bar\lambda$, respectively, with the same fractional order.  
As discussed in \cite{Xu25}, the high order terms in $\wt C_1(t)$ will not affect the perturbation results. So we only consider the matrix
\[
C_1(t)=\lambda I_p+N+tM_{11},
\]
where $M_{11}$ is the (1,1) block of $D_{11}$ defined in (\ref{defd}).

Partition $U_T$ in $U$ as $U_T=\mat{cc}V&V_c\rix$. Following (\ref{strd}),
\[
CV=V(\lambda I_p+N),\quad CV_c=V_c(\lambda I_p+N)^\ast,\quad \mat{c}V^\ast\\V_c^\ast\rix\Delta\mat{cc}V&V_c\rix = \Delta_T=\mat{cc}0&I_p\\I_p&0\rix.
\]
Using the fact that $C$ is $\Delta$-Hermitian,
\[
(\Delta V)^\ast C=(\lambda I_p+N)^\ast (\Delta V)^\ast,\quad
(\Delta V_c)^\ast C=(\lambda I_p+N) (\Delta V_c)^\ast.
\]
Partition
\[
V=\mat{cccc}V_1&V_2&\ldots&V_{m}\rix,\quad 
V_c=\mat{cccc}V_{1}^{(c)}&V_{2}^{(c)}&\ldots&V_{m}^{(c)}\rix,
\]
with the column numbers of $V_j,V_j^{(c)}$ being identical to that of $N_j$, for $j=1,\ldots,m$.
For each $V_j$ and $V_j^{(c)}$, partition further
\eq{defv}
V_j=\mat{cccc}V_{j1}&V_{j2}&\ldots&V_{jj}\rix,\qquad
V_j^{(c)}=\mat{cccc}V_{j1}^{(c)}&V_{j2}^{(c)}&\ldots&V_{jj}^{(c)}\rix,
\en
with the column numbers of $V_{ji},V_{ji}^{(c)}$ being $s_j$ for $i=1,\ldots,j$ and $j=1,\ldots,m$.
Define
\eq{defe}
E_j=\mat{cccc}V_{j1}&V_{j+1,1}&\ldots&V_{m1}\rix,\quad
E_j^{(c)}=\mat{cccc}V_{jj}^{(c)}&V_{j+1,j+1}^{(c)}&\ldots&V_{mm}^{(c)}\rix,
\en
for $j=1,\ldots,m$. Note $E_1$ is a basis matrix of  the right eigenvector space of $C$ corresponding to the eigenvalue $\lambda$ and $E_1^{(c)}$ is a basis matrix of 
the right eigenvector space of $C$ corresponding to $\bar\lambda$.  $\Delta E_1^{(c)}$ 
and $\Delta E_1$ are basis matrices of the left eigenvector spaces of $C$ 
corresponding to the eigenvalue $\lambda$ and $\bar{\lambda}$, respectively. 
From \eqref{defd}, 
\[
D_{11}=\Delta_TU_T^\ast \Delta DU_T.
\]
Then 
\eq{defm11}
M_{11}=V_c^\ast \Delta DV.
\en
Let the matrix 
$B$ in Lemma~\ref{lem1} be $M_{11}$. Then following (\ref{defw}),
$W_m=(V_{mm}^{(c)})^\ast \Delta DV_{m1}$ and
\eq{defhw}
W_k=(E_k^{(c)})^\ast(\Delta D)E_k=\mat{cc}(V_{kk}^{(c)})^\ast\Delta DV_{k1}&
(V_{kk}^{(c)})^\ast\Delta D E_{k+1}\\(E_{k+1}^{(c)})^\ast \Delta D V_{k1}&
W_{k+1}\rix
\en
for $k=1,\ldots,m-1$.
Similarly, for $\wt C_1^\ast(t)$ we consider the matrix
\[
C_1^\ast(t)=(\lambda I_p+N+tM_{11})^\ast=\bar\lambda I_p+N^\ast+tM_{11}^\ast.
\]
Let  
\eq{defpj}
P=\diag(P_1,\ldots,P_m),\quad P_j=\mat{ccc}&&I_{s_j}\\&\adots&\\I_{s_j}&&\rix_{j\times j},
\quad j=1,\ldots,m.
\en
Then
\eq{vh1}
P^\ast C_1^\ast(t) P=\bar\lambda I_p+N+tP^\ast M_{11}^\ast P.
\en
Since $P^\ast M_{11}^\ast P=P^\ast V^\ast (\Delta D)V_cP$, the matrix 
 corresponding to $W_k$ 
in (\ref{defw}) is
\eq{defhww}
\wt W_k=E_k^\ast (\Delta D)E_k^{(c)}=W_k^\ast,
\en
for $k=1,\ldots,m$.
We have the following results.
\begin{theorem}\label{thm1} Let $\Delta$ be Hermitian and invertible.
Suppose both $C$ and $D$ are $\Delta$-Hermitian matrices and $C+tD$ has
the structured decomposition (\ref{nonimag}) when $t$ is sufficiently small. 
Let $W_1,\ldots,W_m$ be defined in (\ref{defhw}) corresponding to the eigenvalue 
$\lambda$ of $C$.
For a fixed $\rho\in\{1,\ldots,m\}$, assume $W_{\rho+1}$ is invertible if $\rho<m$. Let
\[
G_{\rho}=
-W_{\rho+1}^{-1}(E_{\rho+1}^{(c)})^\ast \Delta DV_{\rho1},
\qquad
G^{(c)}_\rho=-W_{\rho+1}^{-\ast}E_{\rho+1}^\ast \Delta DV_{\rho\rho}^{(c)},
 \]
 and let
$\gamma_1^{(\rho)},\ldots,\gamma_{s_\rho}^{(\rho)}$ be the $s_\rho$ eigenvalues of the matrix
\[
S_\rho=\left\{
\begin{array}{ll}
(V_{\rho\rho}^{(c)})^\ast \Delta D(V_{\rho1}+E_{\rho+1} G_{\rho})&\rho<m\\
(V_{mm}^{(c)})^\ast \Delta DV_{m1}&\rho=m,
\end{array}
\right.
\]
where $E_{\rho+1}$, $E_{\rho+1}^{(c)}$, $V_{\rho1}$, $V_{\rho\rho}^{(c)}$
are defined in  (\ref{defe}) and (\ref{defv}), respectively.
\begin{enumerate}
\item[(a)] For each $i\in\{1,\ldots,s_\rho\}$, let $\mu_{i1}^{(\rho)},\ldots,\mu_{i\rho}^{(\rho)}$ be the $\rho$th roots
of $\gamma_i^{(\rho)}$. Then $C+tD$ has $\rho s_\rho$ pairs of eigenvalues $\lambda_{ij}^{(\rho)}(t),
\overline{\lambda_{ij}^{(\rho)}(t)}$ with
\[
\lambda_{ij}^{(\rho)}(t)=\lambda+t^{1/\rho}\mu_{ij}^{(\rho)}+o(t^{1/\rho}),
\quad i=1,\ldots,s_\rho,\quad j=1,\ldots,\rho.
\]
\item[(b)] Suppose all $W_1,\ldots,W_m$ are  invertible. 
Let $\Theta_\rho(S_\rho)$ be defined as in (\ref{deftheta}) with $S_\rho$ given above.
There are matrices
\bstar
\wt V_\rho(t)&=&E_\rho\mat{cc}I_{s_\rho}&0\\G_\rho&0\rix+\sum_{k=1}^mV_k\wt X_{k\rho}(t)+O(t),\\
\wt V_\rho^{(c)}(t)&=&E_\rho^{(c)}\mat{cc}0&I_{s_\rho}\\0&G_\rho^{(c)}\rix
+\sum_{k=1}^mV_k^{(c)}\wt X_{k\rho}^{(c)}(t)+O(t),
\estar
with
\[
\wt X_{k\rho}(t)=
\left\{\begin{array}{ll}
O(\mat{cccc}t^{1-k/\rho}&t^{1-(k-1)/\rho}&\ldots&t^{1-1/\rho}\rix^T)& k<\rho\\
\diag(0_{s_\rho},t^{1/\rho}I_{s_\rho},\ldots,t^{1-1/\rho}I_{s_\rho})&k=\rho\\
O(\mat{ccccccc}t^{1/\rho}&t^{1/\rho}&t^{2/\rho}&\ldots&t&\ldots&t\rix^T)&
k>\rho\end{array}\right.
\]
and
\[
\wt X^{(c)}_{k\rho}(t)=\left\{\begin{array}{ll}
O(\mat{cccc}t^{1-1/\rho}&t^{1-2/\rho}&\ldots&t^{1-k/\rho}\rix^T)&k<\rho\\
\diag(t^{1-1/\rho}I_{s_\rho},\ldots,t^{1/\rho}I_{s_\rho},0_{s_\rho})
+O(\mat{cccc}t&t^{1-1/\rho}&\ldots&t^{1/\rho}\rix^T)&k=\rho\\
O(\mat{ccccccc}
t&\ldots&t&\ldots&t^{2/\rho}&t^{1/\rho}&t^{1/\rho}\rix^T)&k>\rho,
\end{array}\right.
\]
such that the matrix $\wt U_\rho(t)=\mat{cc}\wt V_\rho(t)&\wt V_{\rho}^{(c)}(t)\rix$ satisfies
\eq{jrho}
(C+tD)\wt U_\rho(t)=\wt U_\rho(t)
\mat{cc}\lambda I_{\rho s_\rho}+t^{1/\rho}\wt\Theta_\rho(t)&0\\0&
(\lambda I_{\rho s_\rho}+t^{1/\rho}\wt \Theta_\rho(t))^\ast\rix
\en
and $\wt U^\ast_\rho(t) \Delta \wt U_\rho(t)=t^{1-1/\rho}\Delta_\rho$, where
$\wt \Theta_\rho(t)=\Theta_\rho(S_\rho)+O(t^{1/\rho})$ and
$ \Delta_\rho=\mat{cc}0& I_{\rho s_{\rho}}\\ I_{\rho s_{\rho}}&0\rix$.
\end{enumerate}
\end{theorem}
\proof
Part (a) is directly from Lemma~\ref{lem1}.

For part (b),
by applying Lemma~\ref{lem1} to $\wt C_1(t)$ and using (\ref{nonimag}), one 
has
\[
(C+tD)\wt V_\rho(t)=\wt V_\rho(t) (\lambda I_{\rho s_\rho}+t^{1/\rho}\wt \Theta_\rho(t)),
\]
where $\wt \Theta_\rho(t)=\Theta_\rho(S_\rho)+O(t^{1/\rho})$ and
\eq{deftv}
\wt V_\rho(t) =UZ(t)\mat{c}\wt Z_1(t)\mat{c}X_\rho(t)\\0\rix\\0\rix
\en
with $X_\rho(t)$ given in (\ref{defx}). Since $Z(t)=I_{2n}+O(t)$ and $\wt Z_1(t)=I_{2p}+O(t)$,
\[
\wt V_\rho(t) =U\mat{c}X_\rho(t)\\0\\0\rix+O(t)=VX_\rho(t)+O(t).
\]
The formula for $\wt V_\rho(t)$ is derived by using the block structures of $X_\rho(t)$
and  $V$.

For the block $\wt C_1^\ast(t)$, we apply Lemma~\ref{lem1} to $P^\ast \wt C_1^\ast (t) P$.
Based on (\ref{vh1}), the matrix corresponding to $W_k$ is $\wt W_k=W_k^\ast$ defined in (\ref{defhww}).
Hence, the matrix corresponding to $S_\rho$ is
$
 S_\rho^\ast,
$
and we have
\[
(C+tD)\wh V_\rho^{(c)}(t) =\wh V_\rho^{(c)}(t)(P_\rho^\ast (\bar \lambda I_{\rho s_\rho}+t^{1/\rho}\wt \Theta_\rho^{(c)}(t))P_\rho),
\]
where $\wt \Theta_\rho^{(c)}(t)=\Theta_\rho(S_\rho^\ast)+O(t^{1/\rho})$,
\eq{deftvc}
\wh V_\rho^{(c)}(t)=UZ(t)\mat{c}\wt Z_1(t)\mat{c}0\\P^\ast X^{(c)}_\rho(t) P_\rho\rix\\0\rix
=V_cP^\ast X^{(c)}_\rho(t) P_\rho+O(t),
\en
 $P$ and $P_\rho$ are defined in (\ref{defpj}),
and $X^{(c)}_\rho(t)$ has the same pattern as $X_\rho(t)$ with $G_\rho$ replaced by $G_\rho^{(c)}$. 
Partition
$P^\ast X^{(c)}_\rho(t) P_\rho
=\mat{ccc}(X_{1\rho}^{(c)})^T(t)&\ldots&(X_{m\rho}^{(c)})^T(t)\rix^T$. One
has 
\[
X^{(c)}_{k\rho}(t)=
\left\{\begin{array}{ll}
O(\mat{cccc}t^{1-1/\rho}&t^{1-2/\rho}&\ldots&t^{1-k/\rho}\rix^T)& k<\rho\\
\mat{cc}0&0\\0&I_{s_\rho}\rix
+\diag(t^{1-1/\rho}I_{s_\rho},\ldots,t^{1/\rho}I_{s_\rho},0_{s_\rho})&k=\rho\\
\mat{cc}0&0\\0&G_{k\rho}^{(c)}\rix+O(\mat{ccccccc}
t&\ldots&t&\ldots&t^{2/\rho}&t^{1/\rho}&t^{1/\rho}\rix^T)& k>\rho.
\end{array}\right.
\]

Using the $\Delta$-Hermitian property of $C+tD$, one has
\bstar
(\Delta\wt V_\rho(t))^\ast (C+tD)&=&(\lambda I_{\rho s_\rho}+t^{1/\rho}\wt \Theta_\rho(t))^\ast (\Delta\wt V_\rho(t))^\ast\\
(\Delta\wh V_\rho^{(c)}(t))^\ast (C+tD)&=&(P_\rho^\ast(\bar \lambda I_{\rho s_\rho}+t^{1/\rho} \wt\Theta_\rho^{(c)}( t))P_\rho)^\ast (\Delta\wh V_\rho^{(c)}(t))^\ast.
\estar
This indicates that $\Delta\wt V_\rho(t)$ is a basis matrix of the left invariant subspace of all the eigenvalues
of $C+tD$ that are perturbed from $\bar \lambda$ of order $O(t^{1/\rho})$. 
Similarly, $\Delta\wh V_\rho^{(c)}(t)$ is a basis matrix of the left invariant subspace of all the eigenvalues
of $C+tD$ that are perturbed from $\lambda$ of order $O(t^{1/\rho})$.
Therefore, 
\[
(\Delta\wt V_\rho(t))^\ast \wt V_\rho(t)=(\Delta\wh V_\rho^{(c)}(t))^\ast \wh V_\rho^{(c)}(t)=0,
\]
and $(\Delta\wt V_\rho(t))^\ast \wh V_\rho^{(c)}(t)$ is invertible.
Using (\ref{proz1}), (\ref{deftv}), and (\ref{deftvc}),
\eqn
\nonumber
\wt V_\rho^\ast(t) \Delta\wh V_\rho^{(c)}(t)&=&
\left(UZ(t)\mat{c}\wt Z_1(t)\mat{c}X_\rho(t)\\0\rix\rix\right)^\ast \Delta \left(
UZ(t)\mat{c}\wt Z_1(t)\mat{c}0\\P^\ast X_\rho^{(c)}(t)P_\rho\rix\\0\rix\right)\\
\nonumber
&=&\left(\wt Z_1(t)\mat{c}X_\rho(t)\\0\rix\right)^\ast \Delta_T\left(\wt Z_1(t)\mat{c}0\\P^\ast X_\rho^{(c)}(t)P_\rho\rix
\right)\\
\nonumber
&=&\mat{c}X_\rho(t)\\0\rix^\ast \Delta_T\mat{c}0\\P^\ast X_\rho^{(c)}(t)P_\rho\rix
=X_\rho^\ast(t) P^\ast X_\rho^{(c)}(t)P_\rho
\\\nonumber
&=&\sum_{k=1}^mX_{k\rho}^\ast(t) X_{k\rho}^{(c)}(t)
=t^{1-1/\rho}I_{\rho s_\rho}+O(t)=t^{1-1/\rho}R_\rho(t),
\enn
where $R_\rho(t)=I_{\rho s_\rho}+O(t^{1/\rho})$. 
Let 
\eq{defwhxc}
\wh X_\rho^{(c)}(t)=P^\ast X_\rho^{(c)}(t)P_\rho R_\rho^{-1}(t),\qquad
\wt V_\rho^{(c)}(t)=\wh V_\rho^{(c)}(t)R_\rho^{-1}(t)
=V_c\wh X_\rho^{(c)}(t)+O(t),
\en
and $\wt U_\rho(t)=\mat{cc}\wt V_\rho(t)&\wt V_\rho^{(c)}(t)\rix$.
One has
\[
\wt U_\rho^\ast(t) \Delta \wt U_\rho(t)
=\mat{cc}0&\wt V_\rho^\ast(t) \Delta\wt V_\rho^{(c)}(t)\\
(\wt V_\rho^\ast(t) \Delta\wt V_\rho^{(c)}(t))^\ast&0\rix
=t^{1-1/\rho}\Delta_{\rho}.
\]
Since $R_\rho^{-1}(t)=I_{\rho s_\rho}+O(t^{1/\rho})$, partition 
$\wh X_\rho^{(c)}(t)=\mat{ccc}(\wh X_{1\rho}^{(c)})^T(t)&\ldots&(\wh X_{m\rho}^{(c)})^T(t)\rix^T$ 
conformably 
with $P^\ast X_\rho^{(c)}(t)P_\rho$, one has
\[
\wh X^{(c)}_{k\rho}(t)=X^{(c)}_{k\rho}(t)R_\rho^{-1}(t)=
\left\{\begin{array}{ll}
O(\mat{cccc}t^{1-1/\rho}&t^{1-2/\rho}&\ldots&t^{1-k/\rho}\rix^T)&k<\rho\\
\mat{cc}0&0\\0&I_\rho\rix
+\diag(t^{1-1/\rho}I_{s_\rho},\ldots,t^{1/\rho}I_{s_\rho},0_{s_\rho})&\\
\qquad+O(\mat{cccc}t&t^{1-1/\rho}&\ldots&t^{1/\rho}\rix^T),&k=\rho\\
\mat{cc}0&0\\0&G_{k\rho}^{(c)}\rix+O(\mat{ccccccc}
t&\ldots&t&\ldots&t^{2/\rho}&t^{1/\rho}&t^{1/\rho}\rix^T)& k>\rho,
\end{array}\right.
\]
with which the expression of  $\wt V_\rho^{(c)}(t)$ can be derivd.

Finally,
\[
(C+tD)\wt U_\rho(t)=\wt U_\rho(t)\mat{cc}
\lambda I_{\rho s_\rho}+t^{1/\rho}\wt \Theta_\rho(t)&0\\
0&R_\rho(t) P_\rho^\ast (\bar \lambda I_{\rho s_{\rho}}+t^{1/\rho}\wt\Theta_\rho^{(c)}(t))P_\rho R_\rho^{-1}(t)\rix.
\]
Because
\bstar
&&\wt U_\rho^\ast(t) \Delta(C+tD)\wt U_\rho(t)\\
&&\qquad=\wt U_\rho^\ast (t)\Delta\wt U_\rho(t)
\mat{cc}
\lambda I_{\rho s_\rho}+t^{1/\rho}\wt \Theta_\rho(t)&0\\
0&R_\rho(t) P_\rho^\ast (\bar \lambda I_{\rho s_{\rho}}+t^{1/\rho}\wt\Theta_\rho^{(c)}(t))P_\rho R_\rho^{-1}(t)\rix\\
&&\qquad=t^{1-1/\rho}\Delta_\rho \mat{cc}
\lambda I_{\rho s_\rho}+t^{1/\rho}\wt \Theta_\rho(t)&0\\
0&R_\rho(t) P_\rho^\ast (\bar \lambda I_{\rho s_{\rho}}+t^{1/\rho}\wt\Theta_\rho^{(c)}(t))P_\rho R_\rho^{-1}(t)\rix
\estar
is Hermitian, one has
\[
R_\rho(t) P_\rho^\ast (\bar \lambda I_{\rho s_{\rho}}+t^{1/\rho}\wt\Theta_\rho^{(c)}(t))P_\rho R_\rho^{-1}(t)=(\lambda I_{\rho s_\rho}+t^{1/\rho}\wt \Theta_\rho(t))^\ast.
\]
\eproof

\medskip
The next result is about the perturbations of the invariant subspaces corresponding to 
subsets of  
$\Lambda(\lambda I_{\rho s_\rho}+t^{1/\rho}\wt \Theta_\rho(t))$ and 
$\Lambda((\lambda I_{\rho s_\rho}+t^{1/\rho}\wt \Theta_\rho(t))^\ast)$. 

Suppose $F,F^{(c)}\in\setC^{\rho s_\rho\times \ell}$ are full column rank matrices satisfying
\eq{invth}
\Theta_\rho(S_\rho) F=F\Omega,\quad \Theta_\rho^\ast (S_\rho)F^{(c)}
=F^{(c)}\Omega^\ast,\quad
F^\ast F^{(c)}=I_\ell,
\en
for some $\Omega\in\setC^{\ell\times\ell}$ that satisfies
\eq{proomega}
\Lambda(\Omega)\subseteq\Lambda(\Theta_\rho(S_\rho)),\qquad
\Lambda(\Omega)\cap(\Lambda(\Theta_\rho(S_\rho))\backslash\Lambda(\Omega))=\emptyset.
\en
Due to the block structure of $\Theta_\rho(S_\rho)$, $F$ and $F^{(c)}$ have the block forms
\[
F=\mat{cccc}Q^T&(Q\Omega)^T&\ldots&(Q\Omega^{\rho-1})^T\rix^T,\quad
F^{(c)}=\mat{cccc}\Omega^{\rho-1}(Q^{(c)})^\ast&\ldots&\Omega(Q^{(c)})^\ast&(Q^{(c)})^\ast\rix^\ast,
\]
for some matrices $Q,Q^{(c)}\in\setC^{s_\rho\times \ell}$ satisfying
\[
S_\rho Q=Q\Omega^\rho,\quad S_\rho^\ast Q^{(c)}=Q^{(c)}(\Omega^\ast)^{\rho}.
\]
\begin{theorem}\label{thm2} Under the conditions of Theorem~\ref{thm1}, 
suppose $F$ and $F^{(c)}$ are the full column rank matrices satisfying (\ref{invth}),
and  $\Omega$ satisfies (\ref{proomega}). 
For $t$ sufficiently small, 
there exist matrices $\wt F(t)=F+O(t^{1/\rho})$, $\wt F^{(c)}(t)=F^{(c)}+O(t^{1/\rho})$,
and  $\wt \Omega(t)=\Omega+O(t^{1/\rho})$ satisfying
\[
\wt\Theta_\rho(t)\wt F(t)=\wt F(t)\wt\Omega(t),\quad
\wt\Theta_\rho^\ast(t)\wt F^{(c)}(t)=\wt F^{(c)}(t)\wt\Omega^\ast(t),\quad
\wt F^\ast(t)\wt F^{(c)}(t)=I_\ell,
\]
such that the matrix $\wt U_\Omega(t)=\mat{cc} L_\rho(t)& L_{\rho}^{(c)}(t)\rix$ satisfies
\[
(C+tD)\wt U_\Omega(t)=\wt U_\Omega(t)\mat{cc}\lambda I_\ell+t^{1/\rho}\wt\Omega(t)&0\\0&
(\lambda I_\ell+t^{1/\rho}\wt\Omega(t))^\ast\rix,
\qquad
\wt U^\ast_\Omega(t) \Delta \wt U_\Omega(t)=t^{1-1/\rho}\Delta_\Omega,
\]
where 
$
\Delta_\Omega=\mat{cc}0& I_\ell\\ I_\ell&0\rix,
$
and
\bstar
L_\rho(t)&=&
\wt V_\rho(t)\wt F(t)=E_\rho\mat{cc}I_{s_\rho}\\G_\rho\rix Q+\sum_{k=1}^mV_kL_{k\rho}(t)+O(t),\\
 L_\rho^{(c)}&=&\wt V^{(c)}_\rho(t)\wt F^{(c)}(t)=E_\rho^{(c)}\mat{cc}I_{s_\rho}\\G_\rho^{(c)}\rix Q^{(c)}
+\sum_{k=1}^mV_k^{(c)} L_{k\rho}^{(c)}(t)+O(t)
\estar
with 
\[
 L_{k\rho}(t)=\left\{\begin{array}{ll}
 O(\mat{cccc}t^{1-k/\rho}&t^{1-(k-1)/\rho}&\ldots&t^{1-1/\rho}\rix^T)&k<\rho\\
 \mat{cccc}0&(Q(t^{1/\rho}\Omega))^T&\ldots&(Q(t^{1/\rho}\Omega)^{\rho-1})^T\rix^T
+O(\mat{cccc}t^{1/\rho}&t^{2/\rho}&\ldots&t\rix^T)&k=\rho\\
O(\mat{ccccccc}t^{1/\rho}&t^{1/\rho}&t^{2/\rho}&\ldots&t&\ldots&t\rix^T)&k>\rho
\end{array}\right.
\]
and
\[
L^{(c)}_{k\rho}(t)=\left\{
\begin{array}{ll}O(\mat{cccc}t^{1-1/\rho}&t^{1-2/\rho}&\ldots&t^{1-k/\rho}\rix^T)&k<\rho\\
\mat{cccc}(t^{1/\rho}\Omega)^{\rho-1}(Q^{(c)})^\ast&\ldots&(t^{1/\rho}\Omega)(Q^{(c)})^\ast&0\rix^\ast&\\
\qquad\qquad\qquad\qquad\qquad\qquad+O(\mat{cccc}t&\ldots&t^{2/\rho}&t^{1/\rho}\rix^T)&k=\rho\\
O(\mat{ccccccc}
t&\ldots&t&\ldots&t^{2/\rho}&t^{1/\rho}&t^{1/\rho}\rix^T)&k>\rho.
\end{array}\right.
\]
\end{theorem}
\proof Because 
\[
\wt \Theta_\rho(t)=\Theta_\rho(S_\rho)+O(t^{1/\rho}),
\]
following the standard perturbation theory (\cite{SteS90}), for $t$ sufficiently small, there are 
\[
\wt F(t)=F+O(t^{1/\rho}),\quad \wh F^{(c)}(t)=F^{(c)}+O(t^{1/\rho})
\]
such that 
\[
\wt \Theta_\rho(t)\wt F(t)
=\wt F(t)\wt\Omega(t),\quad
\wt \Theta_\rho^\ast(t))
\wh F^{(c)}(t)=\wh F^{(c)}(t)\wh\Omega^\ast(t),
\]
with 
$
\wt \Omega(t)=\Omega+O(t^{1/\rho}),
$
$\wh\Omega(t)=\Omega+O(t^{1/\rho}).
$
Let 
\[
\wh R(t)=\wt F^\ast(t) \wh F^{(c)}(t)=F^\ast F^{(c)}+O(t^{1/\rho})=I_\ell+O(t^{1/\rho}).
\]
Define 
\[
\wt F^{(c)}(t)=\wh F^{(c)}(t)\wh R^{-1}(t).
\]
Then 
\[
\wt F^\ast(t) \wt F^{(c)}(t)=I_\ell,\quad \wt F^{(c)}(t)=F^{(c)}+O(t^{1/\rho}).
\]
Let $L_\rho(t)=\wt V_\rho(t)\wt F(t)$, $ L_\rho^{(c)}(t)=\wt V_\rho^{(c)}(t)\wt F^{(c)}(t)$,
and $\wt U_\Omega(t)=\mat{cc}L_\rho(t)&L_\rho^{(c)}(t)\rix$. The results can be proved by using (\ref{jrho}).
\eproof

When $\Omega$ is a scalar, we have the following result. 
\begin{corollary}\label{cor1} 
Under the conditions of Theorem~\ref{thm1}, suppose $\omega$ is a simple  eigenvalue of 
$\Theta_\rho(S_\rho)$. The vectors
\[
f=\mat{cccc}q^T&\omega q^T&\ldots&\omega^{\rho-1}q^T\rix^T,\quad
f^{(c)}=\mat{cccc}\omega^{\rho-1}(q^{(c)})^\ast&\ldots&\omega(q^{(c)})^\ast
&(q^{(c)})^\ast\rix^\ast,
\]
satisfy
\[
\Theta_\rho(S_\rho)f=\omega f,\quad
\Theta_\rho^\ast(S_\rho)f^{(c)}=\bar\omega f^{(c)};\quad
S_\rho q=\omega^\rho q,\quad S_\rho^\ast q^{(c)}=\bar\omega^{\rho}q^{(c)},
\]
and
\[
f^\ast f^{(c)}=\rho{\bar\omega}^{\rho-1}q^\ast q^{(c)}=1.
\]
For $t$ sufficiently small, there are vectors $\wt f(t)=f+O(t^{1/\rho})$,
$\wt f^{(c)}(t)=f^{(c)}+O(t^{1/\rho})$, and a  scalar $\wt\omega(t)=\omega+O(t^{1/\rho})$
satisfying
\[
\wt \Theta_\rho(t)\wt f(t)=\wt\omega(t)\wt f(t),\quad
\wt\Theta_\rho^\ast(t)\wt f^{(c)}=\overline{\wt\omega(t)}\wt f^{(c)}(t),\quad
\wt f^\ast(t) \wt f^{(c)}(t)=1,
\]
such that matrix $\wt U_\omega(t)=\mat{cc}l_\rho(t)&l_\rho^{(c)}(t)\rix$ satisfies
\[
(C+tD)\wt U_\omega(t)=\wt U_\omega(t)\mat{cc}\lambda+t^{1/\rho}\wt\omega(t)&0\\0&\overline{\lambda+t^{1/\rho}\wt\omega(t)}\rix,\quad
\wt U_\omega^\ast(t) \Delta \wt U_\omega(t)=t^{1-1/\rho}\mat{cc}0&1\\1&0\rix,
\]
where 
\bstar
l_\rho(t)&=&\wt V_\rho(t)\wt f(t)=
 E_\rho\mat{cc}I_{s_\rho}\\G_\rho\rix q+\sum_{k=1}^mV_kl_{k\rho}(t)+O(t),\\
 l_\rho^{(c)}(t)&=&\wt V_\rho^{(c)}(t)\wt f^{(c)}(t)=E_\rho^{(c)}\mat{cc}I_{s_\rho}\\G_\rho^{(c)}\rix q^{(c)}
+\sum_{k=1}^mV_k^{(c)} l_{k\rho}^{(c)}(t)+O(t)
\estar
with
\[
 l_{k\rho}(t)=
 \left\{\begin{array}{ll}
 O(\mat{cccc}t^{1-k/\rho}&t^{1-(k-1)/\rho}&\ldots&t^{1-1/\rho}\rix^T)& k<\rho\\
\mat{cccc}0&t^{1/\rho}\omega q^T&\ldots&(t^{1/\rho}\omega)^{\rho-1} q^T\rix^T
+O(\mat{cccc}t^{1/\rho}&t^{2/\rho}&\ldots&t\rix^T)&k=\rho\\
O(\mat{ccccccc}t^{1/\rho}&t^{1/\rho}&t^{2/\rho}&\ldots&t&\ldots&t\rix^T)&k>\rho,
\end{array}\right.
\]
\[
l^{(c)}_{k\rho}(t)=\left\{\begin{array}{ll}
O(\mat{cccc}t^{1-1/\rho}&t^{1-2/\rho}&\ldots&t^{1-k/\rho}\rix^T)& k<\rho\\
\mat{cccc}(t^{1/\rho}\omega)^{\rho-1}(q^{(c)})^\ast&\ldots&
t^{1/\rho}\omega (q^{(c)})^\ast&0\rix^\ast+O(\mat{cccc}t&\ldots&t^{2/\rho}&t^{1/\rho}\rix^T)&k=\rho\\
O(\mat{ccccccc}
t&\ldots&t&\ldots&t^{2/\rho}&t^{1/\rho}&t^{1/\rho}\rix^T)& k>\rho.
\end{array}\right.
\]
\end{corollary}
\proof
It is a special case of Theorem~\ref{thm2} with $\ell=1$.
\eproof
\begin{remark}\label{rem32}\rm 
Based on (\ref{nonimag}), the matrix $U\mat{c}I\\Y(t)\rix Z_1(t)\wt Z_1(t)$ is a basis matrix of the invariant subspace of $C+tD$ 
corresponding to the eigenvalues of $\diag(\wt C_1(t),\wt C_1^\ast(t))$ and it satisfies 
\[
\left(U\mat{c}I\\Y(t)\rix Z_1(t)\wt Z_1(t)\right)^\ast \Delta \left(U\mat{c}I\\Y(t)\rix Z_1(t)\wt Z_1(t)\right)=\Delta_T.
\]
If we apply Theorem~\ref{thm1} to every $\rho$ from $1$ to $m$, we are able to assemble another basis matrix $\wt U(t)=\mat{ccc}\wt U_1(t),\ldots,\wt U_m(t)\rix$ corresponding to the eigenvalues
of $\diag(\wt C_1(t),\wt C_1^\ast(t))$ as well. However, these two basis matrices are different.  
The former basis matrix corresponds to $\wt C_1(t)$ and $\wt C_1^\ast(t)$ that are 
perturbed from $\lambda I_p+N$ and $(\lambda I_p+N)^\ast$ with an order $O(t)$
perturbation. On the other hand, $\wt U(t)$ corresponds to a block diagonal matrix with diagonal blocks perturbed from $\lambda I_{\rho s_\rho}$ and $\bar\lambda I_{\rho s_\rho}$
with perturbations of orders $O(t^{1/\rho})$ for $\rho=1,\ldots,m$. Meanwhile $\wt U_\rho(t)=O(1)$,
for $1\le \rho\le m$, and
\[
\wt U^\ast(t) \Delta \wt U(t)=\diag(\Delta_1,t^{1/2}\Delta_2,\ldots,t^{1-1/m}\Delta_m).
\]
The former basis matrix simply
considers all the eigenvalues perturbed from $\lambda$ and $\bar\lambda$ as a single
set. In contrast, $\wt U(t)$ provides information about every individual perturbed eigenvalue,
which is shown in Theorem~\ref{thm2} and Corollary~\ref{cor1}. In fact, the proofs
of Theorem~\ref{thm2} and Corollary~\ref{cor1} reveal that the perturbation formulas for
smaller invariant subspaces and eigenvectors can be 
derived from $\wt U(t)$ with transformations that are essentially constant.  
\end{remark}
\begin{remark}\rm\label{rem33}
In Corollary~\ref{cor1}, $l_\rho(t)$ and $l_\rho^{(c)}(t)$ are right and left eigenvectors of 
$C+tD$ corresponding to the eigenvalue $\lambda+t^{1/\rho}\wt w(t)$. Note that
$l_\rho(t),l_\rho^{(c)}(t)=O(1)$ and
 $l_\rho^\ast(t) l_\rho^{(c)}(t)=t^{1-1/\rho}$. The latter indicates that  the spectral condition
numbers of $C+tD$ corresponding to the eigenvalue $\lambda+t^{1/\rho}\wt \omega(t)$
is $O(t^{-(1-1/\rho)})$. If we consider
$C$ as the matrix perturbed from $C+tD$,  since the spectral condition
number is $O(t^{-(1-1/\rho)})$, the eigenvalue should change with an order $O(t^{1/\rho})$. 
In this sense, Corollary~\ref{cor1} generalizes the standard first order perturbation results for simple eigenvalues (\cite[Sec.7.2.2]{GolV13}).
\end{remark}

\subsection{Real eigenvalue case}\label{sub12}
In this case, we consider (\ref{nrcase}), where
$\wt C_T(t)=C_T+t\wt D_{11}(t)$,
\[
C_T=\alpha I_p+N,\quad \wt D_{11}(t)= D_{11}+O(t)=\Delta_TU_T^\ast \Delta DU_T+O(t),
\]
$\alpha\in\setR$, and $\Delta_T$ is given in (\ref{defcdb}). Recall that $C_T$, $D_{11}$, and $\wt D_{11}(t)$ 
are $\Delta_T$-Hermitian.
Partition
\eq{defu}
U_T=\mat{cccc}U_1&U_2&\ldots&U_m\rix
\en
and
\[
U_j=\mat{cccc}U_{j1}&U_{j2}&\ldots&U_{jj}\rix,\quad j=1,\ldots,m,
\]
conformably with the block columns of $N$ and $N_j$, respectively. Define
\[
E_k=\mat{cccc}U_{k1}&U_{k+1,1}&\ldots&U_{m1}\rix,\quad k=1,\ldots,m.
\]
Then
\[
D_{11}=\Delta_TU_T^\ast \Delta D U_T=\Delta_T\mat{cccc} U_1^\ast \Delta DU_1&U_1^\ast \Delta DU_2&
\ldots&U_1^\ast \Delta DU_m\\
U_2^\ast \Delta DU_1&U_2^\ast \Delta DU_2&
\ldots&U_2^\ast \Delta DU_m\\
\vdots&\vdots&\ddots&\vdots\\
U_m^\ast \Delta DU_1&U_m^\ast \Delta DU_2&
\ldots&U_m^\ast \Delta DU_m\rix.
\]
If we consider $B(t)$ and $B$ in Lemma~\ref{lem1} as  $\wt D_{11}(t)$ and $D_{11}$,
one has
\eq{defrw}
W_k =\mat{ccc}\Sigma_k&&\\&\ddots&\\
&&\Sigma_m\rix\wh W_k,\quad
\wh W_k=E_k^\ast \Delta D E_k=\mat{cc}U_{k1}^\ast \Delta DU_{k1}&U_{k1}^\ast \Delta DE_{k+1}\\
E_{k+1}^\ast \Delta D U_{k1}&\wh W_{k+1}
\rix
\en
for  $k=1,2,\ldots,m-1$, and $\wh W_{m}=U_{m1}^\ast \Delta DU_{m1}$.
Note $\wh W_k$ is Hermitian for any $k$.
\begin{theorem}\label{thm3} Assume both $C$ and $D$ are $\Delta$-Hermitian matrices and $C+tD$ has the decomposition (\ref{nrcase}) when $t$ is sufficiently small, with $C_T=\alpha I_p+N$, $\alpha\in\setR$, and $\Delta_T$ defined in (\ref{defcdb}). 
Let $W_k$ be defined as in (\ref{defrw}) for $k\in\{1,\ldots,m\}$.
For a fixed $\rho\in\{1,\ldots,m\}$ assume $W_{\rho+1}$ is invertible if $\rho<m$. Define
\[
G_{\rho}=-\wh W_{\rho+1}^{-1}E_{\rho+1}^\ast \Delta DU_{\rho1},
\]  
and let
$\gamma_1^{(\rho)},\ldots,\gamma_{s_\rho}^{(\rho)}$ be the $s_\rho$ 
eigenvalues of the matrix
\[S_\rho=\Sigma_\rho\wh S_\rho,\quad\mbox{where}\quad
\wh S_\rho=\wh S_\rho^\ast=\left\{
\begin{array}{ll}
U_{\rho1}^\ast \Delta D(U_{\rho1}+E_{\rho+1}G_{\rho})&\rho<m\\
U_{m1}^\ast\Delta DU_{m1}&\rho=m.
\end{array}
\right.
\]

\begin{enumerate}
\item[(a)] For each $i\in\{1,\ldots,s_\rho\}$, let $\mu_{i1}^{(\rho)},\ldots,\mu_{i\rho}^{(\rho)}$ 
be the $\rho$th roots
of $\gamma_i^{(\rho)}$. Then $C+tD$ has $\rho s_\rho$ eigenvalues
\[
\lambda_{ij}^{(\rho)}(t)=\alpha+t^{1/\rho}\mu_{ij}^{(\rho)}+o(t^{1/\rho}),
\qquad i=1,\ldots,s_\rho,\quad j=1,\ldots,\rho.
\]
\item[(b)] Suppose $W_1,\ldots,W_m$ are all invertible. 
Let $\Theta_\rho(S_\rho)$ be defined as (\ref{deftheta}) with $S_\rho$ given above.
There is a matrix
\[
\wt U_\rho(t)=E_\rho\mat{cc}I_{s_\rho}&0\\G_\rho&0\rix+\sum_{k=1}^mU_k\wt X_{k\rho}(t)+O(t),\\
\]
satisfying
\[
(C+tD)\wt U_\rho(t)=\wt U_\rho(t)
(\alpha I_{\rho s_\rho}+t^{1/\rho}\wt \Theta_\rho(t)),\qquad
\wt \Theta_\rho(t)=\Theta_\rho(S_\rho)+O(t^{1/\rho})
\]
and $\wt U_\rho^\ast(t) \Delta \wt U_\rho(t)= t^{1-1/\rho}\Gamma_{\rho}$,
where
\[
\wt X_{k\rho}(t)=\left\{\begin{array}{ll}
O(\mat{cccc}t^{1-k/\rho}&t^{1-(k-1)/\rho}&\ldots&t^{1-1/\rho}\rix^T)& k<\rho\\
\diag(0_{s_\rho},t^{1/\rho}I_{s_\rho},\ldots,t^{1-1/\rho}I_{s_\rho})
+O(\mat{cccc}t^{1/\rho}&t^{2/\rho}&\ldots&t\rix^T)&k=\rho\\
O(\mat{ccccccc}t^{1/\rho}&t^{1/\rho}&t^{2/\rho}&\ldots&t&\ldots&t\rix^T)& k>\rho,
\end{array}\right.
\]
 $U_1,\ldots,U_m$ are the block columns of $U_T$ given in (\ref{defu}), and
 $\Gamma_\rho$ is defined in Proposition~\ref{pro1}.
\end{enumerate}
\end{theorem}
\proof 
From (\ref{nrcase}) one has
\[
(C+tD)U\mat{c}I\\Y(t)\rix Z_1(t)=U\mat{c}I\\Y(t)\rix Z_1(t)\wt C_T(t),
\]
\[
\left(U\mat{c}I\\Y(t)\rix Z_1(t)\right)^\ast \Delta\left(U\mat{c}I\\Y(t)\rix Z_1(t)\right)= \Delta_T,
\]
where $Y(t)=O(t)$, $Z_1(t)=I_r+O(t^2)$, and $\wt C_T(t)=C_T+t\wt D_{11}(t)$, 
$\wt D_{11}(t)=D_{11}+O(t)$, $D_{11}=\Delta_TU_T^\ast\Delta DU_T$ satisfying
\[
(\Delta_T\wt C_T(t))^\ast=\Delta_T\wt C_T(t),\quad
(\Delta_TD_{11})^\ast =\Delta_TD_{11}.
\]
Applying Lemma~\ref{lem1}(a) to $\wt C_T(t)$ with $B(t)=\wt D_{11}(t)$ and $B=D_{11}$, 
we have part (a). 

For part (b), following Lemma~\ref{lem1}(b) 
one has
\[
(C+tD)\wh U_\rho(t)=\wh U_\rho(t)
(\alpha I_p+t^{1/\rho}\wh \Theta_\rho(t)),
\]
where $\wh \Theta_\rho(t)=\Theta_\rho(S_\rho)+O(t^{1/\rho})$,
\[
\wh U_\rho(t)=U\mat{c}I\\Y(t)\rix Z_1(t)X_\rho(t),
\]
and $X_\rho(t)$ is defined as in (\ref{defx}). 
Using (\ref{nrcase}),
\bstar
\wh U_\rho^\ast(t) \Delta \wh U_\rho(t)
&=&X_\rho ^\ast(t) \Delta_TX_\rho(t)=
X_{\rho\rho}^\ast(t) \Gamma_\rho X_{\rho\rho}(t)
+\sum_{k\ne \rho}X_{k\rho}^\ast (t)\Gamma_kX_{k\rho}(t)\\
&=& t^{1-1/\rho}\Gamma_\rho+O(t)
=t^{1-1/\rho}\wt R_\rho(t),
\estar
where $\wt R_\rho(t)=\Gamma_\rho+O(t^{1/\rho})$. Note $\wt R_\rho(t)$ is Hermitian. As
did in the beginning of the section for (\ref{eqdt}), when $t$ is sufficiently small, one has
$Z_\rho(t)=I+O(t^{1/\rho})$ such that
\[
Z_\rho^\ast(t)\wt R_\rho(t) Z_\rho(t)  =\Gamma_\rho.
\]
Let 
\[
\wt U_\rho(t)=\wh U_\rho(t) Z_\rho(t)=U_TX_\rho(t) Z_\rho(t) +O(t),\quad
\wt \Theta_\rho(t)=Z_\rho^{-1}(t)\wh\Theta_\rho(t)Z_\rho(t)
=\Theta_\rho(S_\rho)+O(t^{1/\rho}).
\]
The results follow.
\eproof

\begin{remark}\label{rem1.5}\rm
From Theorem~\ref{thm3}(b), $\wt \Theta_\rho(t)$ is $\Gamma_\rho$-Hermitian. 
Therefore,  the eigenvalues of 
$\wt\Theta_\rho(t)$, or equivalently 
the eigenvalues $\lambda_{ij}^{(\rho)}(t)$, are either in conjugate pairs or real. Since 
 in this case $S_\rho$ is $\Sigma_\rho$-Hermitian, $\Theta_\rho(S_\rho)$ is $\Gamma_\rho$-Hermitian as well.
When $t$ is sufficiently small, the real eigenvalues of $\wt\Theta_\rho(t)$ are perturbed from 
those real $\mu_{ij}^{(\rho)}$'s, which correspond to the real $\gamma_i^{(\rho)}$'s 
when $\rho$ is odd or real positive $\gamma_i^{(\rho)}$'s when $\rho$ is even. However, because of the relation
$\wt\Theta_\rho(t)=\Theta_\rho(S_\rho)+O(t^{1/\rho})$, not all real $\mu_{ij}^{(\rho)}$'s correspond to the real eigenvalues of $\wt \Theta_\rho(t)$.  
The $O(t^{1/\rho})$ term could introduce nonzero imaginary parts.  
\end{remark}

\medskip
For the perturbations of invariant subspaces with smaller dimensions, because $\Theta_\rho(S_\rho)$
and $\wt \Theta_\rho(t)$ are $\Gamma_\rho$-Hermitian matrices, we will
consider the situation with a full column rank matrix $F\in\setC^{\rho s_\rho\times \ell}$ such that
\eq{condfo}
\Theta_\rho(S_\rho)F=F\Omega,\quad F^\ast \Gamma_\rho F= \Delta_\Omega,
\en
where $\Omega\in\setC^{\ell\times \ell}$ with the eigenvalues either real or in complex conjugate pairs,
$\Lambda(\Omega)\cap(\Lambda(\Theta_\rho(S_\rho))\backslash\Lambda(\Omega))=\emptyset$,
and  $\Delta_\Omega$ is Hermitian and block diagonal as $\wh \Delta$
defined in Proposition~\ref{pro1}. 
Due to the block structure of $\Theta_\rho(S_\rho)$,
it still has
\eq{formf}
F=\mat{cccc}Q^T&(Q\Omega)^T&\ldots&(Q\Omega^{\rho-1})^T\rix^T
\en
for some matrix $Q\in\setC^{s_\rho\times \ell}$ satisfying $S_\rho Q=Q\Omega^\rho$.
Using the same method as we did in the beginning of the section, we can show
the existence of 
\eq{ordfo}
\wt F(t)=F+O(t^{1/\rho}),\qquad \wt\Omega(t)=\Omega+O(t^{1/\rho})
\en
satisfying
\eq{pertf}
\wt \Theta_\rho(t)\wt F(t)=\wt F(t)\wt\Omega(t),\quad \wt F^\ast(t)\Gamma_\rho \wt F(t)= \Delta_\Omega,
\en

%
\begin{theorem}\label{thm4} Under the conditions of Theorem~\ref{thm3}, 
suppose $F$ is a full column rank matrix of the form (\ref{formf}) that satisfies (\ref{condfo}).
 Also, suppose $\Lambda(\Omega)\subseteq\Lambda(\Theta_\rho(S_\rho))$ with
the eigenvalues in $\Lambda(\Omega)$ being either real 
or in complex conjugate pairs and $\Lambda(\Omega)\cap(\Lambda(\Theta_\rho(S_\rho))\backslash\Lambda(\Omega))=\emptyset$. Then for $t$ sufficiently small, 
there exist matrices $\wt F(t)$ and $\wt \Omega(t)$
satisfying (\ref{ordfo}) and (\ref{pertf}),
such that the matrix
\[
\wt U_\Omega(t)=\wt U_\rho(t)\wt F(t)=
E_\rho\mat{cc}I_{s_\rho}\\G_\rho\rix Q+\sum_{k=1}^mU_kL_{k\rho}(t)+O(t)
\]
with
\[
 L_{k\rho}(t)=\left\{\begin{array}{ll}
 O(\mat{cccc}t^{1-k/\rho}&t^{1-(k-1)/\rho}&\ldots&t^{1-1/\rho}\rix^T)& k<\rho\\
\mat{cccc}0&(Q(t^{1/\rho}\Omega))^T&\ldots&(Q(t^{1/\rho}\Omega)^{\rho-1})^T\rix^T+O(\mat{cccc}t^{1/\rho}&t^{2/\rho}&\ldots&t\rix^T)&k=\rho\\
O(\mat{ccccccc}t^{1/\rho}&t^{1/\rho}&t^{2/\rho}&\ldots&t&\ldots&t\rix^T)& k>\rho
\end{array}\right.
\]
satisfies
\[
(C+tD)\wt U_\Omega(t)=\wt U_\Omega(t)(\alpha I_\ell+t^{1/\rho}\wt \Omega(t)),\qquad
\wt U^\ast_\Omega(t) \Delta \wt U_\Omega(t)=t^{1-1/\rho}\Delta_\Omega.
\]
\end{theorem}
\proof
The results can be derived by setting
$
\wt U_\Omega(t)=\wt U_\rho(t)\wt F(t),
$
where $\wt U_\rho(t)$ is given in Theorem~\ref{thm3}, and using (\ref{ordfo})
and (\ref{pertf}).
\eproof

\medskip
The following results are the special cases when either $\Omega$ is a  simple real eigenvalue
of $\Theta_\rho(S_\rho)$, or $\Omega=\mat{cc}\mu&0\\0&\bar\mu\rix$, 
where $\mu,\bar\mu$ are nonreal simple eigenvalues of $\Theta_\rho(S_\rho)$. 
\begin{corollary}\label{cor2} 
\begin{enumerate}
\item[(a)] Suppose $\omega$ is a simple real eigenvalue of 
$\Theta_\rho(S_\rho)$. There is a vector
\eq{def2f}
f=\mat{cccc}q^T&\omega q^T&\ldots&\omega^{\rho-1}q^T\rix^T
\en
that satisfies 
\eq{pro2f}
\Theta_\rho(S_\rho)f=\omega f,\quad S_\rho q=\omega^\rho q, \quad f^\ast \Gamma_\rho f
=\rho\omega^{\rho-1}q^\ast\Sigma_\rho q=\sigma,\quad \sigma=\pm 1.
\en
When
 $t$ is sufficiently small,
 there exist a vector $\wt f(t)=f+O(t^{1/\rho})$ and a real scalar $\wt\omega(t)=\omega+O(t^{1/\rho})$
 satisfying
 \[
 \wt \Theta_\rho(t)\wt f(t)=\wt\omega(t)\wt f(t),\qquad
 \wt f^\ast(t) \Gamma_\rho \wt f(t)=\sigma,
 \]
 such that the vector
 \eq{formwtu}
 \wt u_\omega(t) = \wt U_\rho(t) \wt f(t)=E_\rho\mat{cc}I_{s_\rho}\\G_\rho\rix q+\sum_{k=1}^mU_kl_{k\rho}(t)+O(t),
 \en
 with
\[
 l_{k\rho}(t)=\left\{\begin{array}{ll}
 O(\mat{cccc}t^{1-k/\rho}&t^{1-(k-1)/\rho}&\ldots&t^{1-1/\rho}\rix^T)& k<\rho\\
 \mat{cccc}0&t^{1/\rho}\omega q^T&\ldots&(t^{1/\rho}\omega)^{\rho-1} q^T\rix^T+O(\mat{cccc}t^{1/\rho}&t^{2/\rho}&\ldots&t\rix^T)&k=\rho\\
O(\mat{ccccccc}t^{1/\rho}&t^{1/\rho}&t^{2/\rho}&\ldots&t&\ldots&t\rix^T)& k>\rho
\end{array}\right.
\]
satisfies
\[
(C+tD)\wt u_\omega(t)=(\alpha I+t^{1/\rho}\wt \omega(t))\wt u_\omega(t),
\qquad 
\wt u_\omega^\ast (t)\Delta \wt u_\omega(t)=t^{1-1/\rho}\sigma.
\]
\item[(b)] Suppose
\[
\Omega=\mat{cc}\mu&0\\0&\bar\mu\rix,
\]
where $\mu,\bar\mu$ are nonreal simple eigenvalues of $\Theta_\rho(S_\rho)$. There are
vectors
\[
f=\mat{cccc}q^T&\mu q^T&\ldots&\mu^{\rho-1}q^T\rix^T,\quad
f^{(c)}=\mat{cccc}(q^{(c)})^\ast&\mu(q^{(c)})^\ast&\ldots&\mu^{\rho-1}(q^{(c)})^\ast\rix^\ast.
\]
that satisfy
\[
\Theta_\rho(S_\rho)f=\mu f,\quad
\Theta_\rho(S_\rho)f^{(c)}=\bar\mu f^{(c)},\quad
S_\rho q=\mu^\rho q,\quad S_\rho q^{(c)}=\bar\mu^\rho q^{(c)},
\]
and
\[(f^{(c)})^\ast \Gamma_\rho f^{(c)}=f^\ast \Gamma_\rho f=0,\qquad
f^\ast \Gamma_\rho f^{(c)}=\rho(\bar\mu)^{\rho-1}q^\ast\Sigma_\rho q^{(c)}=1.
\]
For $t$ sufficiently small, there are vectors $\wt f(t)=f+O(t^{1/\rho})$, $\wt f^{(c)}(t)=f^{(c)}+O(t^{1/\rho})$,
and a nonreal scalar $\wt\mu(t)=\mu+O(t^{1/\rho})$ satisfying
\[
\wt\Theta_\rho(t)\wt f(t)=\wt\mu(t)\wt f(t),\quad
\wt\Theta_\rho(t)\wt f^{(c)}(t)=\overline{\wt\mu(t)}\wt f^{(c)}(t),
\]
\[
\wt f^\ast(t)\Gamma_\rho \wt f(t)=(\wt f^{(c)}(t))^\ast\Gamma_\rho \wt f^{(c)}(t)=0,\qquad
\wt f^\ast(t)\Gamma_\rho \wt f^{(c)}(t)=1.
\]
Then the matrix $\wt U_\Omega(t)=\mat{cc}l_\rho(t)&l_\rho^{(c)}(t)\rix$ satisfies
\[
(C+tD)\wt U_\Omega(t)=\wt U_\Omega(t)\mat{cc}\alpha+t^{1/\rho}\wt\mu(t)&0\\0&
\overline{\alpha+t^{1/\rho}\wt\mu(t)}\rix,\quad
\wt U_\Omega^\ast(t) \Delta \wt U_\Omega(t)=t^{1-1/\rho}\mat{cc}0&1\\1&0\rix,
\]
where 
\bstar
l_\rho(t)&=&\wt U_\rho(t) \wt f(t)=E_\rho\mat{cc}I_{s_\rho}\\G_\rho\rix q+\sum_{k=1}^mU_kl_{k\rho}(t)+O(t),\\
 l_\rho^{(c)}(t)&=&\wt U_\rho(t)\wt f^{(c)}(t)=E_\rho\mat{cc}I_{s_\rho}\\G_\rho\rix q^{(c)}
+\sum_{k=1}^mU_k l_{k\rho}^{(c)}(t)+O(t)
\estar
with
\[
 l_{k\rho}(t)=\left\{\begin{array}{ll}
 O(\mat{cccc}t^{1-k/\rho}&t^{1-(k-1)/\rho}&\ldots&t^{1-1/\rho}\rix^T)& k<\rho\\
 \mat{cccc}0&t^{1/\rho}\mu q^T&\ldots&(t^{1/\rho}\mu)^{\rho-1} q^T\rix^T+O(\mat{cccc}t^{1/\rho}&t^{2/\rho}&\ldots&t\rix^T)&k=\rho\\
O(\mat{ccccccc}t^{1/\rho}&t^{1/\rho}&t^{2/\rho}&\ldots&t&\ldots&t\rix^T)& k>\rho
\end{array}\right.
\]
and
\[
l_{k\rho}^{(c)}(t)=\left\{\begin{array}{ll}
O(\mat{cccc}t^{1-k/\rho}&t^{1-(k-1)/\rho}&\ldots&t^{1-1/\rho}\rix^T)&k<\rho\\
\mat{cccc}0& t^{1/\rho}\mu (q^{(c)})^\ast&\ldots&
(t^{1/\rho}\mu)^{\rho-1} (q^{(c)})^\ast\rix^\ast&\\
\qquad\qquad\qquad\qquad+O(\mat{cccc}t^{1/\rho}&t^{2/\rho}&\ldots&t\rix^T)
&k=\rho\\
O(\mat{ccccccc}t^{1/\rho}&t^{1/\rho}&t^{2/\rho}&\ldots&t&\ldots&t\rix^T)& k>\rho.
\end{array}\right.
\]
\end{enumerate}
\end{corollary}
\proof
(a) Since $\omega$ is
a simple real eigenvalue of the $\Gamma_\rho$-Hermitian matrix $\Theta_\rho(S_\rho)$, 
following Proposition~\ref{pro1} there is $f\ne 0$ of the form (\ref{def2f}) that satisfies (\ref{pro2f}). Then there is 
\[
\wh f(t) =f+O(t^{1/\rho}),
\]
such that $\wt \Theta_\rho(t)\wh f(t)=\wt \omega(t) \wh f(t)$ with
$\wt \omega(t)=\omega+O(t^{1/\rho})$, and for $\wt U_\rho(t)$ given in Theorem~\ref{thm3},
\[ 
\wh f^\ast(t) \wt U_\rho(t)\Delta \wt U_\rho(t) \wh f(t)=\wh f^\ast(t) (t^{1-1/\rho}\Gamma_\rho) \wh f(t)
=t^{1-1/\rho}\wt \sigma,\qquad \wt \sigma=\sigma+O(t^{1/\rho}).
\]
Note $\wt \sigma$ is real.
When $t$ is sufficiently small, there is $r=1+O(t^{1/\rho})\in \setR$ such that $\wt\sigma r^2=\sigma$.
By setting $\wt f(t)=\wh f(t) r=f+O(t^{1/\rho})$, and $\wt u_\omega(t)=\wt U_\omega(t) \wt f(t)$, one has
\[
(C+tD)\wt u_\omega(t)=(\alpha+t^{1/\rho}\wt \omega(t)) \wt u_\omega(t),\qquad
\wt u_\omega^\ast (t)\Delta \wt u_\omega(t)=t^{1-1/\rho}\sigma,
\]
and $\wt u_\omega(t)$ has the formula as (\ref{formwtu}). Since $C+tD$ is $\Delta$-Hermitian, 
the above relations show that $\wt \omega(t)$ must be real.

(b)  Let $f\ne 0$ satisfy
\[
\Theta_\rho(S_\rho) f =\mu f.
\]
Because $\Theta_\rho(S_\rho)$ is $\Gamma_\rho$-Hermitian and $\mu$ is simple, following
Proposition~\ref{pro1}, there is  $f^{(c)}$ such that both $f,f^{(c)}$ have the given forms, and
$F=\mat{cc}f&f^{(c)}\rix$ satisfies
\[
\Theta_\rho(S_\rho)F=F\mat{cc}\mu&0\\0&\bar\mu\rix,\quad
F^\ast \Gamma_\rho F=\mat{cc}0&1\\1&0\rix.
\]
Similarly, one can show the existence of $\wt f(t)$, 
$\wt f^{(c)}(t)$, and $\wt \mu(t)$. 
The results can be derived by setting $\wt U_\Omega(t)=\mat{cc}l_\rho(t)&l^{(c)}_\rho(t)\rix
=\wt U_\rho(t) \wt F(t)$, where $\wt F(t)=\mat{cc}\wt f(t)&\wt f^{(c)}(t)\rix=F+O(t^{1/\rho})$ and
satisfies
\[
\wt \Theta_\rho(t) \wt F(t)=\wt F(t)\mat{cc}\wt \mu(t)&0\\0&\overline{\wt\mu(t)}\rix,\quad
\wt \mu(t)=\mu+O(t^{1/\rho}),\quad \wt F^\ast(t) \Gamma_\rho \wt F(t)=\mat{cc}0&1\\1&0\rix.
\]
\eproof
%
\section{Invariant subspace perturbations for Hamiltonian matrices}\label{sec4}
We now consider the invariant subspace perturbations for Hamiltonian matrices.
We assume a Hamiltonian matrix $H$ is perturbed to $H+tT$,
where $t$ is sufficiently small and $T$ is another Hamiltonian matrix.
 We will also focus on the following two
cases:
\begin{enumerate}
\item[(a)] a single defective nonimaginary eigenvalue $\lambda$ (paired with $-\bar\lambda$)
\item[(b)] a single defective imaginary eigenvalue $\imath \alpha$, where $\alpha\in\setR$.
\end{enumerate}
Note $H$ is Hamiltonian if and only if $-\imath H$ is
$-\imath J_n$-Hermitian.  
If  we set $C=-\imath H$ , $D=-\imath T$, and $\Delta=-\imath J_n$, we still have (\ref{strd}), (\ref{defcda}),
(\ref{defcdb}). The only difference is that $\Delta$ is congruent to $\diag(I_n,-I_n)$ in this situation.
Replacing the nonreal eigenvalue $\lambda$ in (\ref{defcda}) by $-\imath\lambda$ with a 
nonimaginary $\lambda$,
we have
\[
C_T=-\imath \wh H_T, \quad \Delta_T=U_T^\ast (-\imath J_n)U_T,
\]
where 
\[
\begin{array}{lll}
\wh H_T=\mat{cc}\lambda I_p+\imath N&0\\0&-(\lambda I_p+\imath N)^\ast\rix,&
\Delta_T=\mat{cc}0&I_p\\I_p&0\rix,
&\mbox{case (a)}\\
\wh H_T=\imath(\alpha I+N),&\Delta_T=\diag(\Gamma_1,\ldots,\Gamma_m),&\mbox{case (b)}
\end{array}
\]
and $\Gamma_1,\ldots,\Gamma_m$ are defined in (\ref{defcdb}).
 In terms of $H$ and $J_n$, one has
\eq{hamfac}
HU=U\mat{cc}\wh H_T&0\\0&\wh H_C\rix,\quad
U^\ast J_n U=\imath\mat{cc}\Delta_T&0\\0&\Delta_C\rix,\quad
U=\mat{cc}U_T&U_C\rix,
\en
where $U$ is the same as (\ref{strd}),
 $\wh H_T$ is defined above, and $\wh H_C=\imath C_C$.
\subsection{Nonimaginary eigenvalue case}
Based on the equivalence relations described above, we will use the results developed in Subsection~\ref{sub11}.
In this case, $U_T=\mat{cc}V&V_c\rix$ satisfies
\[
HU_T=U_T\wh H_T=U_T\mat{cc}\lambda I_p+\imath N&0\\0&
-(\lambda I_p+\imath N)^\ast\rix,\qquad
U_T^\ast J_nU_T=\imath\Delta_T=\imath \mat{cc}0&I_p\\I_p&0\rix.
\]
Define  the unitary block diagonal matrix
\eq{defwhp}
\wh P=\diag(\wh P_1,\ldots,\wh P_m),
\en
where
\[
\wh P_j=\diag(I_{s_j}, (-\imath)I_{s_j},\ldots,(-\imath)^{j-1} I_{s_j}),\quad j=1,\ldots,m.
\]
One has
\[
\wh P^\ast (\imath N)\wh P= N.
\]
Hence 
\[
H_T:=\mat{cc}\wh P&0\\0&-\imath\wh P\rix^\ast
\wh H_T
\mat{cc}\wh P&0\\0&-\imath\wh P\rix=\mat{cc}\lambda I_p+N&0\\0&-(\lambda I_p+N)^\ast\rix,
\]
and for $\Phi_T=\mat{cc}\Xi&\Xi_c\rix$ with 
\[
\Xi=V\wh P, \qquad \Xi_c=V_c(-\imath \wh P),
\]
 one has
\[
H\Phi_T=\Phi_TH_T,\qquad \Phi_T^\ast J_n\Phi_T=J_p.
\]
Since $-\imath J_n$ and $\Delta_T$ are congruent to $\diag(I_n,-I_n)$ and $\diag(I_p,-I_p)$,
respectively, $\Delta_C$ in (\ref{hamfac}) is congruent to $\diag(I_{n-p},-I_{n-p})$. Therefore,
there is a $\wh P_C$ such that 
\eq{defpc}
\wh P_C^\ast(\imath\Delta_C)\wh P_C=J_{n-p}.
\en
Let $\Phi_C=U_C\wh P_C$, $\Phi=\mat{cc}\Phi_T&\Phi_C\rix$
and $H_C=\wh P_C^\ast \wh H_C\wh P_C$. Then (\ref{hamfac}) becomes
\eq{fach}
H\Phi =\Phi\mat{cc}H_T&0\\0&H_C\rix,\quad
\Phi^\ast J_n\Phi =\mat{cc}J_p&0\\0&J_{n-p}\rix.
\en
Partition
\eq{defxibc}
\Xi=\mat{ccc}\Xi_1&\ldots&\Xi_m\rix,\qquad
\Xi_c=\mat{ccc}\Xi_1^{(c)}&\ldots&\Xi_m^{(c)}\rix,
\en
\[
\Xi_j=\mat{ccc}\Xi_{j1}&\ldots&\Xi_{jj}\rix,\quad
\Xi_j^{(c)}=\mat{ccc}\Xi_{j1}^{(c)}&\ldots&\Xi_{jj}^{(c)}\rix,\quad
j=1,\ldots,m,
\]
consistent with the block columns of $N$ and $N_j$, respectively.
Set
\eq{defup}
\Upsilon_k=\mat{ccc}\Xi_{k1}&\ldots&\Xi_{m1}\rix,\qquad
\Upsilon_k^{(c)}=\mat{ccc}\Xi_{kk}^{(c)}&\ldots&\Xi_{mm}^{(c)}\rix.
\en
Tracing back to the matrix $C+tD=-\imath(H+tT)$, 
the block $M_{11}$ defined in (\ref{defm11}) is
\[
M_{11}=V_c^\ast (-J_n T)V=(\imath\wh P)\Xi_c^\ast (J_nT)\Xi\wh P^\ast
\]
and
\[
W_k=-\mat{ccc}(-\imath)^kI_{s_k}&&\\&\ddots&\\&&(-\imath)^m I_{s_m}\rix \Psi_k,
\]
where
\eq{defpsi}
\Psi_k=(\Upsilon_k^{(c)})^\ast(J_nT)\Upsilon_k=
\mat{ccc}(\Xi_{kk}^{(c)})^\ast (J_nT)\Xi_{k1}&\ldots&(\Xi_{kk}^{(c)})^\ast(J_nT)\Xi_{m1}\\
\vdots&\ddots&\vdots\\
(\Xi_{mm}^{(c)})^\ast(J_nT)\Xi_{k1}&\ldots&(\Xi_{mm}^{(c)})^\ast(J_nT)\Xi_{m1}\rix.
\en
Then
\eq{relggs}
G_\rho=\wh G_\rho,\quad G_\rho^{(c)}=\mat{ccc}-\imath I_{s_\rho}&&\\&\ddots&\\&&
(-\imath)^{m-\rho}I_{s_m}\rix
\wh G^{(c)}_\rho,\quad S_\rho =-(-\imath)^\rho \wh S_\rho
\en
with
\eqn
\nonumber 
&&\wh G_{\rho}=- \Psi_{\rho+1}^{-1}(\Upsilon_{\rho+1}^{(c)})^\ast J_nT\Xi_{\rho1},\quad
\wh G^{(c)}_\rho=
-\Psi_{\rho+1}^{-\ast}\Upsilon_{\rho+1}^\ast J_nT\Xi_{\rho\rho}^{(c)},\\
\label{hggs}
&&\wh S_\rho=\left\{
\begin{array}{ll}(\Xi_{\rho\rho}^{(c)})^\ast J_nT(\Xi_{\rho1}+\Upsilon_{\rho+1}\wh G_\rho)
&\rho<m\\
(\Xi_{mm}^{(c)})^\ast J_nT\Xi_{m1}&\rho =m.
\end{array}\right.
\enn
Based on the results in Subsection~\ref{sub11}, 
we have the following results for Hamiltonian matrices.
%
\begin{theorem}\label{thm5} Let $H$ and $T$ be Hamiltonian matrices
 and $H$ satisfy (\ref{fach})
with $\lambda$ being nonimaginary,  $\Phi=\mat{cc}\Phi_T&\Phi_C\rix$, $\Phi_T=\mat{cc}\Xi&\Xi_c\rix$.
Suppose $t$ is sufficiently small. 
Let $\Psi_1,\ldots,\Psi_m$ be defined  in (\ref{defpsi}).
For a fixed $\rho\in\{1,\ldots,m\}$ assume $\Psi_{\rho+1}$ is invertible if $\rho<m$,
and let
$\gamma_1^{(\rho)},\ldots,\gamma_{s_\rho}^{(\rho)}$ be the $s_\rho$ eigenvalues of the matrix
$-\wh S_\rho$, where $\wh S_\rho$ is defined in (\ref{hggs}).
\begin{enumerate}
\item[(a)] For each $i\in\{1,\ldots,s_\rho\}$, let $\mu_{i1}^{(\rho)},\ldots,\mu_{i\rho}^{(\rho)}$ be the $\rho$th roots
of $\gamma_i^{(\rho)}$. Then $H+tT$ has $\rho s_\rho$ pairs of eigenvalues $\lambda_{ij}^{(\rho)}(t),
-\overline{\lambda_{ij}^{(\rho)}(t)}$ with
\[
\lambda_{ij}^{(\rho)}(t)=\lambda+t^{1/\rho}\mu_{ij}^{(\rho)}+o(t^{1/\rho}),
\quad i=1,\ldots,s_\rho,\quad j=1,\ldots,\rho.
\]

\item[(b)] Suppose all $\Psi_1,\ldots, \Psi_m$ are invertible. 
Let $\Theta_\rho(-\wh S_\rho)$ be defined as (\ref{deftheta}) with $S_\rho$ replaced by
$-\wh S_\rho$.
Then for the matrices
\bstar
\wt \Xi_\rho(t)&=&\Upsilon_\rho\mat{cc}I_{s_\rho}&0\\\wh G_\rho&0\rix+\sum_{k=1}^m\Xi_k\wt X_{k\rho}(t)+O(t),\\
\wt \Xi_\rho^{(c)}(t)&=&\Upsilon_\rho^{(c)}\mat{cc}0&I_{s_\rho}\\0&\wh G_\rho^{(c)}\rix
+\sum_{k=1}^m\Xi_k^{(c)}\wt X_{k\rho}^{(c)}(t)+O(t),
\estar
with
\[
\wt X_{k\rho}(t)=\left\{\begin{array}{ll}
O(\mat{cccc}t^{1-k/\rho}&t^{1-(k-1)/\rho}&\ldots&t^{1-1/\rho}\rix^T)& k<\rho\\
\diag(0_{s_\rho},t^{1/\rho}I_{s_\rho},\ldots,t^{1-1/\rho}I_{s_\rho})&k=\rho
\\
O(\mat{ccccccc}t^{1/\rho}&t^{1/\rho}&t^{2/\rho}&\ldots&t&\ldots&t\rix^T)& k>\rho
\end{array}\right.
\]
and
\[
\wt X^{(c)}_{k\rho}(t)=\left\{\begin{array}{ll}
O(\mat{cccc}t^{1-1/\rho}&t^{1-2/\rho}&\ldots&t^{1-k/\rho}\rix^T)& k<\rho\\
\diag(t^{1-1/\rho}I_{s_\rho},\ldots,t^{1/\rho}I_{s_\rho},0_{s_\rho})
+O(\mat{cccc}t&t^{1-1/\rho}&\ldots&t^{1/\rho}\rix^T)&k=\rho\\
O(\mat{ccccccc}
t&\ldots&t&\ldots&t^{2/\rho}&t^{1/\rho}&t^{1/\rho}\rix^T)& k>\rho,
\end{array}\right.
\]
where $\Xi_k, \Xi_k^{(c)}$ and $\Upsilon_\rho,\Upsilon_\rho^{(c)}$  are defined in (\ref{defxibc}) and (\ref{defup}), respectively,
the matrix $\wt \Phi_\rho(t)=\mat{cc}\wt \Xi_\rho(t)&\wt \Xi_{\rho}^{(c)}(t)\rix$ satisfies
\[
(H+tT)\wt \Phi_\rho(t)=\wt \Phi_\rho(t)
\mat{cc}\lambda I_{\rho s_\rho}+t^{1/\rho}\wt\Theta_\rho^{(h)}(t)&0\\0&
-(\lambda I_{\rho s_\rho}+t^{1/\rho}\wt \Theta_\rho^{(h)}(t))^\ast\rix
\]
and $\wt \Phi^\ast_\rho(t) J_n \wt \Phi_\rho(t)=t^{1-1/\rho}J_{\rho s_\rho}$, where
$\wt \Theta_\rho^{(h)}(t)=\Theta_\rho(-\wh S_\rho)+O(t^{1/\rho})$.
\end{enumerate}
\end{theorem}
\proof
In Theorem~\ref{thm1}, set $C=-\imath H$, $D=-\imath T$, $\Delta =-\imath J_n$,
and replace $\lambda$ by $-\imath\lambda$ and $U$ by $\Phi=U\diag(\wh P,-\imath \wh P, \wh P_C)$,
where $\wh P$ and $\wh P_C$ are defined in (\ref{defwhp}) and (\ref{defpc}).
For $t$ sufficiently small, from (\ref{jrho}) one has 
\[
(H+tT)\wt\Phi_\rho (t)=\wt\Phi_\rho(t)\mat{cc}\lambda I_{\rho s_\rho}+t^{1/\rho}\wt \Theta_\rho^{(h)}(t)&0\\0&
-(\lambda I_{\rho s_\rho}+t^{1/\rho}\wt \Theta_\rho^{(h)}(t))^\ast\rix
\]
with (by using (\ref{relggs}))
\[
\wt\Theta_\rho^{(h)}(t):=\wh P_\rho^\ast(\imath \wt \Theta_\rho(S_\rho))\wh P_\rho
=\Theta_\rho(\imath^\rho S_\rho)+O(t^{1/\rho}) = \Theta_\rho(-\wh S_\rho)+O(t^{1/\rho}) 
\]
and
\[
\wt\Phi_\rho^\ast(t) J_n \wt\Phi_\rho (t)=t^{1-1/\rho}J_p,
\]
where 
\[
\wt \Phi_\rho(t)=\wt U_\rho(t)\mat{cc}\wh P_\rho&0\\0&-\imath\wh P_\rho\rix
=\mat{cc}\wt\Xi_\rho(t)&\wt\Xi_\rho^{(c)}(t)\rix
\] 
and
\bstar
\wt \Xi_\rho(t) &=&\wt V_\rho(t)\wh P_\rho=VX_\rho(t)\wh P_\rho+O(t)
=\Xi \wh P^\ast X_\rho(t)\wh P_\rho+O(t),\\
\wt\Xi_\rho^{(c)}(t)&=&\wt V^{(c)}_\rho(t)(-\imath \wh P_\rho)
=V_c \wh X_\rho^{(c)}(t)(-\imath\wh P_\rho)+O(t)=\Xi_c \wh P^\ast \wh X_\rho^{(c)}(t)\wh P_\rho+O(t),
\estar
with $X_\rho(t)$ defined in (\ref{defx}) and $\wh X_\rho^{(c)}(t)$  in (\ref{defwhxc}). 
The results can be proven with the above relations
and the block forms of $X_\rho(t)$ and $\wh X_\rho^{(c)}(t)$.
\eproof
\begin{theorem}\label{thm6} Under the conditions of Theorem~\ref{thm5}, 
let $F,F^{(c)}\in\setC^{\rho s_\rho\times \ell}$ be full column rank matrices satisfying 
\bstar
F&=&\mat{cccc}Q^T&(Q\Omega)^T&\ldots&(Q\Omega^{\rho-1})^T\rix^T,\quad
-\wh S_\rho Q=Q\Omega^\rho,\\
F^{(c)}&=&\mat{cccc}\Omega^{\rho-1}(Q^{(c)})^\ast&\ldots&\Omega(Q^{(c)})^\ast&(Q^{(c)})^\ast\rix^\ast,
\quad -\wh S_\rho^\ast Q^{(c)}=Q^{(c)}(\Omega^\ast)^{\rho}
\estar
and
\[
\Theta_\rho(-\wh S_\rho)F=F\Omega,\quad 
\Theta_\rho^\ast (-\wh S_\rho)F^{(c)}=F^{(c)}\Omega^\ast,\quad F^\ast F^{(c)}=
I_\ell
\] 
with the condition $\Lambda(\Omega)\cap(\Lambda(\Theta_\rho(-\wh S_\rho))\backslash\Lambda(\Omega))=\emptyset$.
Then for $t$ sufficiently small, there exist matrices $\wt F(t)=F+O(t^{1/\rho})$,
 $\wt F^{(c)}(t)=F^{(c)}+O(t^{1/\rho})$,
and $\wt \Omega(t)=\Omega+O(t^{1/\rho})$ satisfying
\[
\wt \Theta_\rho^{(h)}(t)\wt F(t)=\wt F(t)\wt\Omega(t),\quad
(\wt \Theta_\rho^{(h)}(t))^\ast\wt F^{(c)}(t)=\wt F^{(c)}(t)\wt\Omega^\ast(t),\qquad
\wt F^\ast \wt F^{(c)}=I_\ell,
\]
such that  $\wt \Phi_\Omega(t)=\mat{cc}L_\rho(t)&L_\rho^{(c)}(t)\rix$ satisfies
\[
(H+tT)\wt \Phi_\Omega(t)=\wt \Phi_\Omega(t)\mat{cc}\lambda I_\ell+t^{1/\rho}\wt\Omega(t)&0\\0&
-(\lambda I_\ell+t^{1/\rho}\wt\Omega(t))^\ast\rix,\qquad
\wt \Phi^\ast_\Omega(t) J_n \wt \Phi_\Omega(t)=t^{1-1/\rho}J_\ell,
\]
where
\bstar
L_\rho(t)&=&
\wt \Xi_\rho(t)\wt F(t)=\Upsilon_\rho\mat{cc}I_{s_\rho}\\\wh G_\rho\rix Q+\sum_{k=1}^m\Xi_kL_{k\rho}(t)+O(t),\\
 L_\rho^{(c)}(t)&=&\wt\Xi_\rho^{(c)}(t)\wt F^{(c)}(t)=\Upsilon_\rho^{(c)}\mat{cc}I_{s_\rho}\\\wh G_\rho^{(c)}\rix Q^{(c)}
+\sum_{k=1}^m\Xi_k^{(c)} L_{k\rho}^{(c)}(t)+O(t)
\estar
with
\[
 L_{k\rho}(t)=\left\{\begin{array}{ll}
 O(\mat{cccc}t^{1-k/\rho}&t^{1-(k-1)/\rho}&\ldots&t^{1-1/\rho}\rix^T)&k<\rho\\
\mat{cccc}0&(Q(t^{1/\rho}\Omega))^T&\ldots&(Q(t^{1/\rho}\Omega)^{\rho-1})^T\rix^T&\\
\qquad\qquad\qquad\qquad\qquad+O(\mat{cccc}t^{1/\rho}&t^{2/\rho}&\ldots&t\rix^T)&k=\rho\\
O(\mat{ccccccc}t^{1/\rho}&t^{1/\rho}&t^{2/\rho}&\ldots&t&\ldots&t\rix^T)& k>\rho
\end{array}\right.
\]
and
\[
L^{(c)}_{k\rho}(t)=\left\{\begin{array}{ll}
O(\mat{cccc}t^{1-1/\rho}&t^{1-2/\rho}&\ldots&t^{1-k/\rho}\rix^T)& k<\rho\\
\mat{cccc}(t^{1/\rho}\Omega)^{\rho-1}(Q^{(c)})^\ast&\ldots&(t^{1/\rho}\Omega)(Q^{(c)})^\ast&0\rix^\ast&\\
\qquad\qquad\qquad\qquad\qquad
+O(\mat{cccc}t&\ldots&t^{2/\rho}&t^{1/\rho}\rix^T)&k=\rho\\
O(\mat{ccccccc}
t&\ldots&t&\ldots&t^{2/\rho}&t^{1/\rho}&t^{1/\rho}\rix^T)& k>\rho.
\end{array}\right.
\]
\end{theorem}
\proof
The proof is essentially the same as that for Theorem~\ref{thm2}.
\eproof
%
%
\begin{corollary}\label{cor3} Suppose $\omega$ is a simple eigenvalue of 
$\Theta_\rho(-\wh S_\rho)$, and the vectors
\[
f=\mat{cccc}q^T&\omega q^T&\ldots&\omega^{\rho-1}q^T\rix^T,\quad
f^{(c)}=\mat{cccc}\omega^{\rho-1}(q^{(c)})^\ast&\ldots&\omega(q^{(c)})^\ast&(q^{(c)})^\ast\rix^\ast
\]
satisfy
\[
\Theta_\rho(-\wh S_\rho)f=\omega f,\quad \Theta_\rho^\ast(-\wh S_\rho) f^{(c)}=\bar\omega f^{(c)},
\quad
f^\ast f^{(c)}=\rho\bar\omega^{\rho-1}q^\ast q^{(c)}=1.
\]
For $t$ sufficiently small, there are vectors
\[
\wt f(t)=f+O(t^{1/\rho}),\quad \wt f^{(c)}(t)=f^{(c)}+O(t^{1/\rho}),
\]
and a scalar $\wt \omega(t)=\omega+O(t^{1/\rho})$ satisfying
\[
\wt \Theta_\rho^{(h)}(t)\wt f(t)=\wt\omega(t)\wt f(t),\qquad
(\wt\Theta_\rho^{(h)}(t))^\ast\wt f^{(c)}(t)=\overline{\wt \omega(t)}\wt f^{(c)}(t),\quad
\wt f^\ast(t)\wt f^{(c)}(t)=1,
\]
such that  $\wt \Phi_\omega(t)=\mat{cc}l_\rho(t)&l_\rho^{(c)}(t)\rix$ satisfies
\[
(H+tT)\wt \Phi_\omega(t)=\wt \Phi_\omega(t)\mat{cc}\lambda+t^{1/\rho}\wt\omega(t)&0\\0&-\overline{(\lambda+t^{1/\rho}\wt\omega(t))}\rix,\qquad
\wt \Phi_\omega^\ast(t) J_n \wt \Phi_\omega(t)=t^{1-1/\rho}J_1,
\]
and
\bstar
l_\rho(t)&=&\wt \Xi_\rho(t)\wt f(t)
=\Upsilon_\rho\mat{cc}I_{s_\rho}\\\wh G_\rho\rix q+\sum_{k=1}^m\Xi_kl_{k\rho}(t)+O(t),\\
 l_\rho^{(c)}(t)&=&\wt \Xi_\rho^{(c)}(t)\wt f^{(c)}(t)
 =\Upsilon_\rho^{(c)}\mat{cc}I_{s_\rho}\\\wh G_\rho^{(c)}\rix q^{(c)}
+\sum_{k=1}^m\Xi_k^{(c)} l_{k\rho}^{(c)}(t)+O(t)
\estar
with
\[
 l_{k\rho}(t)=\left\{\begin{array}{ll}
 O(\mat{cccc}t^{1-k/\rho}&t^{1-(k-1)/\rho}&\ldots&t^{1-1/\rho}\rix^T)&k<\rho\\
\mat{cccc}0&t^{1/\rho}\omega q^T&\ldots&(t^{1/\rho}\omega)^{\rho-1} q^T\rix^T+O(\mat{cccc}t^{1/\rho}&t^{2/\rho}&\ldots&t\rix^T)&k=\rho\\
O(\mat{ccccccc}t^{1/\rho}&t^{1/\rho}&t^{2/\rho}&\ldots&t&\ldots&t\rix^T)& k>\rho,
\end{array}\right.
\]
\[
l^{(c)}_{k\rho}(t)=\left\{\begin{array}{ll}
O(\mat{cccc}t^{1-1/\rho}&t^{1-2/\rho}&\ldots&t^{1-k/\rho}\rix^T)&k<\rho\\
\mat{cccc}(t^{1/\rho}\omega)^{\rho-1}(q^{(c)})^\ast&\ldots&
t^{1/\rho}\omega (q^{(c)})^\ast&0\rix^\ast+O(\mat{cccc}t&\ldots&t^{2/\rho}&t^{1/\rho}\rix^T)&k=\rho\\
O(\mat{ccccccc}
t&\ldots&t&\ldots&t^{2/\rho}&t^{1/\rho}&t^{1/\rho}\rix^T)& k>\rho.
\end{array}\right.
\]
\end{corollary}
\subsection{Purely imaginary eigenvalue case} 
The results will be derived directly from those in Subsection~\ref{sub12} with $C=-\imath H$,
$D=-\imath T$ and $\Delta=-\imath J_n$. For consistency, we replace $U=\mat{cc}U_T&U_C\rix$ 
with $\Phi=\mat{cc}\Phi_T&\Phi_C\rix$. Partition
\eq{defphit}
\Phi_T=\mat{cccc}\Phi_1&\Phi_2&\ldots&\Phi_m\rix,
\en
\[
\Phi_j=\mat{cccc}\Phi_{j1}&\Phi_{j2}&\ldots&\Phi_{jj}\rix,\quad j=1,\ldots,m,
\]
that are consistent with the block columns of $N$ and $N_j$, respectively, and set
\eq{defiup}
\Upsilon_\rho=\mat{ccc}\Phi_{\rho1}&\ldots&\Phi_{m1}\rix
\en
for $\rho\in\{1,\ldots,m\}$. We still have \eqref{hamfac}, and in this case
$\wh H_T=\imath(\alpha I_p+N)$ with $\alpha\in\setR$, and $\Delta_T$ is defined in (\ref{defcdb}).
With the new notations, 
the matrix $W_k$  corresponding to $C$ defined in (\ref{defrw}) becomes
\eq{defhwp}
W_k=-\mat{ccc}\Sigma_k&&\\&\ddots&\\&&\Sigma_m\rix\Psi_k,\qquad \Psi_k=\Upsilon_k^\ast (J_nT)\Upsilon_k= \Psi_k^\ast.
\en
Hence
\eqn
\nonumber
G_\rho&=&-\Psi_{\rho+1}^{-1}\Upsilon_{\rho+1}^\ast(J_nT)\Phi_{\rho1},\\
\label{defhgs}&&\\
\nonumber
S_\rho&=&-\Sigma_\rho\wh S_\rho,\qquad
\wh S_\rho=\wh S_\rho^\ast=\left\{
\begin{array}{ll}
\Phi_{\rho1}^\ast(J_nT)(\Phi_{\rho1}+\Upsilon_{\rho+1}G_\rho)&\rho<m\\
\Phi_{m1}^\ast (J_nT)\Phi_{m1}&\rho=m.
\end{array}
\right.
\enn

\begin{theorem}\label{thm7} Suppose both $H$ and $T$ are Hamiltonian matrices.
For a given $\rho\in\{1,\ldots,m\}$, 
assume $\Psi_{\rho+1}$ defined in (\ref{defhwp}) is invertible if $\rho<m$.
Let $G_\rho$ and $S_\rho$ be defined in (\ref{defhgs})
and 
$\gamma_1^{(\rho)},\ldots,\gamma_{s_\rho}^{(\rho)}$ be the $s_\rho$ 
eigenvalues of the matrix $S_\rho$. 

\begin{enumerate}
\item[(a)] For each $i\in\{1,\ldots,s_\rho\}$, let $\mu_{i1}^{(\rho)},\ldots,\mu_{i\rho}^{(\rho)}$ 
be the $\rho$th roots
of $\gamma_i^{(\rho)}$. Then for $t$ sufficiently small, $H+tT$ has $\rho s_\rho$ eigenvalues
\[
\lambda_{ij}^{(\rho)}(t)=\imath(\alpha+t^{1/\rho}\mu_{ij}^{(\rho)})+o(t^{1/\rho}),
\qquad i=1,\ldots,s_\rho,\quad j=1,\ldots,\rho.
\]
\item[(b)] Suppose all $\Psi_1,\ldots,\Psi_m$ are invertible. 
Let $\Theta_\rho(S_\rho)$ be defined as (\ref{deftheta}) with $S_\rho$ given in (\ref{defhgs}).
Then for $t$ sufficiently small, there is a matrix
\[
\wt \Phi_\rho(t)=\Upsilon_\rho\mat{cc}I_{s_\rho}&0\\G_\rho&0\rix+\sum_{k=1}^m\Phi_k\wt X_{k\rho}(t)+O(t),\\
\]
with $G_\rho$ given in (\ref{defhgs}) and
\[
\wt X_{k\rho}(t)=\left\{\begin{array}{ll}
O(\mat{cccc}t^{1-k/\rho}&t^{1-(k-1)/\rho}&\ldots&t^{1-1/\rho}\rix^T)& k<\rho\\
\diag(0_{s_\rho},t^{1/\rho}I_{s_\rho},\ldots,t^{1-1/\rho}I_{s_\rho})
+O(\mat{cccc}t^{1/\rho}&t^{2/\rho}&\ldots&t\rix^T)&k=\rho\\
O(\mat{ccccccc}t^{1/\rho}&t^{1/\rho}&t^{2/\rho}&\ldots&t&\ldots&t\rix^T)& k>\rho
\end{array}\right.
\]
where $\Phi_1,\ldots,\Phi_m$ are the block columns of $\Phi_T$ defined in (\ref{defphit}) and 
$\Upsilon_\rho$ is defined in (\ref{defiup}), such that
\[
(H+tT)\wt \Phi_\rho(t)=\wt \Phi_\rho(t)
(\imath(\alpha I_{\rho s_\rho}+t^{1/\rho}\wt \Theta_\rho^{(h)}(t))),\qquad
\wt \Theta_\rho^{(h)}(t)=\Theta_\rho(S_\rho)+O(t^{1/\rho})
\]
and $\wt\Phi_\rho^\ast(t) J_n\wt \Phi_\rho(t)= t^{1-1/\rho}(\imath\Gamma_{\rho})$.
\end{enumerate}
\end{theorem}
\proof
It is directly from Theorem~\ref{thm3}.
\eproof
%
\begin{theorem}\label{thm8} Under the conditions of Theorem~\ref{thm7}, 
suppose $F\in\setC^{\rho s_\rho\times \ell}$ is a full column rank matrix of the form (\ref{formf}) that satisfies (\ref{condfo}),
where $\Lambda(\Omega)\subseteq \Lambda(\Theta_\rho(S_\rho))$. 
Assume the eigenvalues of
$\Omega$ are either real 
or in complex conjugate pairs, and $\Lambda(\Omega)\cap(\Lambda(\Theta_\rho(S_\rho))\backslash\Lambda(\Omega))=\emptyset$.
Then   for $t$ sufficiently small,
there is a matrix of the form
\[
\wt \Phi_\Omega(t)=
\Upsilon_\rho\mat{cc}I_{s_\rho}\\G_\rho\rix Q+\sum_{k=1}^m\Phi_kL_{k\rho}(t)+O(t)
\]
with
\[
 L_{k\rho}(t)=\left\{\begin{array}{ll}
 O(\mat{cccc}t^{1-k/\rho}&t^{1-(k-1)/\rho}&\ldots&t^{1-1/\rho}\rix^T)& k<\rho\\
\mat{cccc}0&(Q(t^{1/\rho}\Omega))^T&\ldots&(Q(t^{1/\rho}\Omega)^{\rho-1})^T\rix^T+O(\mat{cccc}t^{1/\rho}&t^{2/\rho}&\ldots&t\rix^T)&k=\rho\\
O(\mat{ccccccc}t^{1/\rho}&t^{1/\rho}&t^{2/\rho}&\ldots&t&\ldots&t\rix^T)&k>\rho
\end{array}\right.
\]
such that 
\[
(H+tT)\wt \Phi_\Omega(t)=\wt \Phi_\Omega(t)(\imath(\alpha I+t^{1/\rho}\wt \Omega(t))),\qquad
\wt \Omega(t)=\Omega+O(t^{1/\rho}),
\]
and $\wt \Phi^\ast_\Omega (t)J_n \wt \Phi_\Omega(t)=t^{1-1/\rho}(\imath\Delta_\Omega)$ with
$ \Delta_\Omega$  given in (\ref{condfo}).
\end{theorem}
\proof
It is directly from Theorem~\ref{thm4}.
\eproof

%
%
\begin{corollary}\label{cor4} 
\begin{enumerate}
\item[(a)] Suppose $\omega$ is a simple real eigenvalue of 
$\Theta_\rho(S_\rho)$. There is a vector
\[
f=\mat{cccc}q^T&\omega q^T&\ldots&\omega^{\rho-1}q^T\rix^T
\]
that satisfies
\[
\Theta_\rho(S_\rho)f=\omega f,\quad S_\rho q=\omega^\rho q,\quad
f^\ast \Gamma_\rho f=\rho\omega^{\rho-1} q^\ast\Sigma_\rho q=\sigma,\quad
\sigma=\pm 1.
\]
 When
 $t$ is sufficiently small,
 there is a vector  $\wt f(t)=f+O(t^{1/\rho})$ and a real scalar $\wt\omega(t)=\omega+O(t^{1/\rho})$
 satisfying
 \[
 \wt \Theta_\rho^{(h)}(t)\wt f(t)=\wt\omega(t)\wt f(t),\qquad \wt f^\ast(t) \Gamma_\rho\wt f(t)=\sigma,
 \]
and the vector
\[
 \wt \phi_\omega(t) =\wt\Phi_\rho(t) \wt f(t)= \Upsilon_\rho\mat{cc}I_{s_\rho}\\G_\rho\rix q+\sum_{k=1}^m\Phi_kl_{k\rho}(t)+O(t),
 \]
 with
\[
 l_{k\rho}(t)=\left\{\begin{array}{ll}
 O(\mat{cccc}t^{1-k/\rho}&t^{1-(k-1)/\rho}&\ldots&t^{1-1/\rho}\rix^T)&k<\rho\\
\mat{cccc}0&t^{1/\rho}\omega q^T&\ldots&(t^{1/\rho}\omega)^{\rho-1} q^T\rix^T+O(\mat{cccc}t^{1/\rho}&t^{2/\rho}&\ldots&t\rix^T)&k=\rho\\
O(\mat{ccccccc}t^{1/\rho}&t^{1/\rho}&t^{2/\rho}&\ldots&t&\ldots&t\rix^T)&k>\rho
\end{array}\right.
\]
satisfies
\[
(H+tT)\wt \phi_\omega(t)=\imath(\alpha I+t^{1/\rho}\wt \omega(t))\wt \phi_\omega(t),\quad
\wt \phi_\omega^\ast(t) J_n \wt \phi_\omega(t)=t^{1-1/\rho}\imath\sigma.
\]
\item[(b)] Suppose
\[
\Omega=\mat{cc}\mu&0\\0&\bar\mu\rix,
\]
where $\mu,\bar\mu$ are nonreal simple eigenvalues of $\Theta_\rho(S_\rho)$.
There are vectors
\[
f=\mat{cccc}q^T&\mu q^T&\ldots&\mu^{\rho-1}q^T\rix^T,\quad
f^{(c)}=\mat{cccc}(q^{(c)})^\ast&\mu (q^{(c)})^\ast&\ldots&\mu^{\rho-1}(q^{(c)})^\ast\rix^\ast
\]
that satisfy
\[
\Theta_\rho(S_\rho)f=\mu f,\quad S_\rho q=\mu^\rho q,\quad
\Theta_\rho(S_\rho) f^{(c)}=\bar\mu f^{(c)}, \quad S_\rho q^{(c)}=\bar\mu^\rho q^{(c)},
\]
and
\[
(f^{(c)})^\ast \Gamma_\rho f^{(c)}=f^\ast\Gamma_\rho f=0,\qquad
f^\ast \Gamma_\rho f^{(c)}=\rho\bar u^{\rho-1}q^\ast\Sigma_\rho q^{(c)}=1.
\]
When $t$ is sufficiently small, there are vectors $\wt f(t)=f+O(t^{(1/\rho})$, 
$\wt f^{(c)}(t)=f^{(c)}+O(t^{(1/\rho})$, and
a nonreal scalar $\wt\mu(t)=\mu+O(t^{(1/\rho})$ satisfying
\[
\wt\Theta_\rho^{(h)}(t)\wt f(t)=\wt\mu(t)\wt f(t),\qquad
\wt\Theta_\rho^{(h)}(t)\wt f^{(c)}(t)=\overline{\wt \mu(t)}\wt f^{(c)}(t),
\]
\[ \wt f^\ast(t) \Gamma_\rho \wt f(t)=(\wt f^{(c)}(t))^\ast \Gamma_\rho \wt f^{(c)}(t)=0,\quad
\wt f^\ast(t)\Gamma_\rho \wt f^{(c)}(t)=1,
\]
such that 
$\wt \Phi_\Omega(t)=\mat{cc}l_\rho(t)&l_\rho^{(c)}(t)\rix$ satisfies
\bstar
(H+tT)\wt \Phi_\Omega(t)&=&\wt \Phi_\Omega(t)\mat{cc}\imath(\alpha+t^{1/\rho}\wt\mu(t))&0\\0&
-\,\overline{[\imath(\alpha+t^{1/\rho}\wt\mu(t))]}
\rix,\\
\wt \Phi_\Omega^\ast (t)J_n\wt \Phi_\Omega(t)&=&\imath t^{1-1/\rho}\mat{cc}0&1\\1&0\rix,
\estar
where 
\bstar
l_\rho(t)&=&
\wt\Phi_\rho (t)\wt f(t)=\Upsilon_\rho\mat{cc}I_{s_\rho}\\G_\rho\rix q+\sum_{k=1}^m\Phi_kl_{k\rho}(t)+O(t),\\
 l_\rho^{(c)}(t)&=&\wt\Phi_\rho(t) \wt f^{(c)}(t)=\Upsilon_\rho\mat{cc}I_{s_\rho}\\G_\rho\rix q^{(c)}
+\sum_{k=1}^m\Phi_k l_{k\rho}^{(c)}(t)+O(t)
\estar
with
\[
 l_{k\rho}(t)=\left\{\begin{array}{ll}
 O(\mat{cccc}t^{1-k/\rho}&t^{1-(k-1)/\rho}&\ldots&t^{1-1/\rho}\rix^T)& k<\rho\\
 \mat{cccc}0&t^{1/\rho}\mu q^T&\ldots&(t^{1/\rho}\mu)^{\rho-1} q^T\rix^T+O(\mat{cccc}t^{1/\rho}&t^{2/\rho}&\ldots&t\rix^T)&k=\rho\\
O(\mat{ccccccc}t^{1/\rho}&t^{1/\rho}&t^{2/\rho}&\ldots&t&\ldots&t\rix^T)& k>\rho
\end{array}\right.
\]
and
\[
l_{k\rho}^{(c)}(t)=\left\{\begin{array}{ll}
O(\mat{cccc}t^{1-1/\rho}&t^{1-2/\rho}&\ldots&t^{1-k/\rho}\rix^T)& k<\rho\\
\mat{cccc}(t^{1/\rho}\mu)^{\rho-1}(q^{(c)})^\ast&\ldots&
t^{1/\rho}\omega (q^{(c)})^\ast&0\rix^\ast&\\
\qquad\qquad\qquad\qquad+O(\mat{cccc}t&\ldots&t^{2/\rho}&t^{1/\rho}\rix^T)&k=\rho\\
O(\mat{ccccccc}
t&\ldots&t&\ldots&t^{2/\rho}&t^{1/\rho}&t^{1/\rho}\rix^T)& k>\rho.
\end{array}\right.
\]
\end{enumerate}
\end{corollary}
\subsection{A special case}
For Hamiltonian matrices an important question is under what conditions the eigenvalues
will stay on the imaginary axis under (small) Hamiltonian perturbations and how to remove them 
from the imaginary axis with Hamiltonian perturbations. Based on the results established above (and Remark~\ref{rem1.5}), when $t$ is sufficiently small,
it is obvious that necessary conditions for
$H+tT$ to have imaginary eigenvalues include
\begin{enumerate}
\item $H$ has imaginary eigenvalues
\item For an imaginary eigenvalue of $H$, the corresponding $S_\rho=-\Sigma_\rho\wh S_\rho$ defined in (\ref{defhgs})  should have real roots when $\rho$ is odd, or  positive real roots when $\rho$ is even.
\end{enumerate}
If the second condition is not satisfied, none of the $\rho$th roots of the eigenvalues
of $S_\rho$ will be real. On the other hand, existence of the real roots of $S_\rho$ is not sufficient
for $H+tT$ to have imaginary eigenvalues. It depends on  $\Sigma_\rho$
and $\wh S_\rho$, and also the higher fractional order term in $\wt\Theta_\rho^{(h)}(t)$. 
 The second condition 
is sufficient if $S_\rho$ has a  semi-simple real eigenvalue (positive when $\rho$ is even)
 with the same corresponding sign signatures. 
  A special case is $t_\rho\ne s_\rho-t_\rho$ in $\Sigma_\rho$ (\cite{MehX99,Tho91}). 
An extreme case is $\Sigma_\rho=I_{s_\rho}$ or $-I_{s_\rho}$.

In \cite{MehX24} a special case was considered with
\[
 H =\wh U \mat{cc|cc}\imath\alpha I_p+N&0&G&0\\0&F_2&0&G_2\\\hline
0&0&-(\imath\alpha I_p+N)^\ast&0\\0&R_2&0&-F_2^\ast\rix \wh U^{-1},
\quad T=J_nK,
\]
where $\alpha\in\setR$, $\wh U$ is  symplectic ($\wh U^\ast J_n\wh U=J_n$), 
$H_C:=\mat{cc}F_2&G_2\\R_2&-F_2^\ast\rix$ is Hamiltonian and 
$\imath\alpha\notin\Lambda(H_C)$. $K$ is positive 
semidefinite, $N$ is defined in (\ref{defn}), and
\[
G=\diag(G_1,\ldots,G_m),\quad G_k=\diag(0,\ldots,0,I_{s_k}),\quad k=1,\ldots,m.
\]
After a  simple two-sided block permutation, we have a nonsingular matrix 
$U=\mat{cc}U_T&U_C\rix$ 
such that
\[
HU=U\mat{cc}H_T&0\\0&H_C\rix,\quad U^\ast J_nU=\mat{cc}\wh J_p&0\\0&J_{n-p}\rix,
\]
where
\[
H_T=\diag(\imath \alpha I_{2s_1}+\wh H_{2},\imath \alpha I_{4s_2}+\wh H_4,\ldots,
\imath \alpha I_{2ms_m}+\wh H_{2m}),
\quad \wh J_p=\diag(J_{s_1},J_{2s_2},\ldots,J_{ms_m})
\]
with
\[
\wh H_{2k}=\mat{cc}N_k&G_k\\0&-N_k^\ast\rix,\quad k=1,\ldots,m.
\]
Let $\wh P_{2k}$ be defined in (\ref{defwhp}) and
\[
\Pi_T=\diag(\Pi_1,\ldots,\Pi_m),
\]
where
\bstar
\Pi_k&=&\mat{cc}I_{ks_k}&0\\0&\begin{array}{ccc}&&(-I_{s_k})^{k-1}\\
&\adots&\\I_{s_k}&&\end{array}\rix \wh P_{2k}^\ast\\
&=&\diag\left(I_{s_k},\imath I_{s_k},\ldots, \imath^{k-1}I_{s_k},\mat{ccc}
&&(-1)^{k-1}\imath^{2k-1}I_{s_k}\\&\adots&\\
\imath^k I_{s_k}&&\rix\right).
\estar
Then
\[
\Pi_T^\ast H_T\Pi_T=\imath(\alpha I+ \wh N), \quad
\Pi_T^\ast \wh J_p\Pi_T=\imath\wh\Delta_T,
\]
where
\[
\wh N=\diag(\wh N_{2},\wh N_{4},\ldots,\wh N_{2m}),\quad
\wh\Delta_T=\imath\diag( \Gamma_2, \Gamma_4,\ldots, \Gamma_{2m})
\]
with
\[
\wh N_{2k}=\mat{cccc}0&I_{s_k}&&\\&\ddots&\ddots&\\&&0&I_{s_k}\\&&&0\rix_{2k\times 2k},\quad
\Gamma_{2k}=\mat{ccc}&&I_{s_k}\\&\adots&\\I_{s_k}&&\rix_{2k\times 2k}.
\]
Note in this special case the sizes of the Jordan blocks of $\imath\alpha$ are all even.
Let 
\[
\Phi=\mat{cc}U_T&U_C\rix\diag(\Pi_T,I)=\mat{cc}\Phi_T&\Phi_C\rix,
\] 
and partition
\[
\Phi_T=\mat{cccc}\Phi_2&\Phi_4&\ldots&\Phi_{2m}\rix,\quad \Phi_{2j}=\mat{cccc}
\Phi_{2j,1}&\Phi_{2j,2}&\ldots&\Phi_{2j,2j}\rix
\]
conformably with the block columns of $\wh N$ and $\wh N_{2k}$. Let
\[
\Upsilon_{2k}=\mat{cccc}\Phi_{2k,1}&\Phi_{2(k+1),1}&\ldots&\Phi_{2m,1}\rix.
\]
Since in this case, $\Sigma_{2k}=I_{s_k}$ for $k=1,\ldots,m$
 and $T=J_nK$, following (\ref{defhwp})
 one has
\[
W_{2k}=-\Upsilon_{2k}^\ast (J_n(J_nK))\Upsilon_{2k}^\ast
=\Upsilon_{2k}^\ast K\Upsilon_{2k}.
\]
Since $K\ge 0$, the nonsingularity requirement becomes $W_{2k}>0$ for $k=1,\ldots.m$. 
Then $S_{2\rho}>0$ for any
$\rho\in\{1,\ldots.m\}$.  Let $\gamma_j>0$
be a semisimple eigenvalue of $S_{2\rho}$ with algebraic (geometric) multiplicity $\ell$.
Then $\Theta_{2\rho}(S_{2\rho})$ has exactly two real semisimple eigenvalues corresponding 
to $\gamma_j$ with the same algebraic multiplicity $\ell$, denoted by
$\mu_{j,\pm}=\pm\gamma_j^{1/(2\rho)}$. Consider $\mu_{j,+}$. Let $Q$ be a basis matrix
of the eigenvector space of $S_{2\rho}$ corresponding to $\gamma_j$, i.e.,
\[
S_{2\rho}Q=\gamma_j Q.
\]
Then 
\[
F=\mat{cccc}Q^T&\mu_{j,+}Q^T&\ldots&\mu_{j,+}^{2\rho-1}Q^T\rix^T
\]
has full column rank and satisfies
\[
\Theta_{2\rho}(S_{2\rho})F=F\Omega,\qquad \Omega=\mu_{j,+}I_\ell.
\]
This shows that $\mu_{j,+}$ is semisimple and $F$ is a basis matrix of the corresponding
eigenvector space of $\Theta_{2\rho}(S_{2\rho})$. Because
\[
\Delta_\Omega:= F^\ast \Gamma_{2\rho}F
=2\rho \mu_{j,+}^{2\rho-1}Q^\ast Q>0.
\]
One can always choose $Q$ such that $\Delta_\Omega=I_{\ell}$.
When $t>0$ is sufficiently small, following Theorem~\ref{thm8}
there are $\wt \Phi_\Omega(t)$ and 
\[
\wt\Omega(t)=\Omega+O(t^{1/(2\rho)})=\mu_{j,+}I_\ell+O(t^{1/(2\rho)})
\]
such that
\[
(H+tJ_nK)\wt \Phi_\Omega(t)=\wt\Phi_\Omega(t) (\imath(\alpha I_\ell+t^{1/(2\rho)} \wt \Omega(t))),\quad
\wt \Phi_\Omega^\ast (t)J_n\wt\Phi_\Omega(t) =\imath t^{1-1/(2\rho)}I_\ell.
\]
Notice 
\[
\wt \Phi_\Omega^\ast (t)J_n(H+tJ_nK)\wt \Phi_\Omega(t)
=-t^{1-1/(2\rho)}(\alpha I_\ell+t^{1/(2\rho)}\wt\Omega(t)).
\]
Since the  left-hand side is Hermitian, so is  $\wt \Omega(t)$. This implies that
$H+tJ_nK$ has $\ell$ (semi)simple purely imaginary  eigenvalues of the form
$
\imath(\alpha+t^{1/(2\rho)}\mu_{j,+}+o(t^{1/(2\rho)}))
$
and the corresponding structure inertia indices are all $1$. In the same way, one can show 
that $H+tJ_nK$ has
$\ell$ purely imaginary eigenvalues of the form
$\imath(\alpha+t^{1/(2\rho)}\mu_{j,-}+o(t^{1/(2\rho)}))$ and the corresponding structure inertia 
indices are all $-1$. The same result was derived in \cite{MehX24}, which now can be
considered as a special case of Theorem~\ref{thm8}.

%
\section{Conclusions}\label{con}\label{sec5} 
By using the recently developed perturbation theory in \cite{Xu25}, we have derived structured 
invariant subspace perturbation
results for both the $\Delta$-Hermitian matrices and Hamiltonian matrices corresponding to defective
eigenvalues. The results exhibit the fractional orders and symmetric structures of the perturbed 
eigenvalues and invariant subspaces, and provide the information that are not contained in
the standard perturbation theory. The perturbation results only apply to the case where all the 
blocks $W_1,\ldots,W_m$ are invertible. Currently, it is not clear how to extend the results to general case.

\medskip
\noindent
{\bf Acknowledgement}
The author thanks the anonymous referee and the handling editor for their comments and suggestions.
\bibliographystyle{amsplain}

\end{document}